\documentclass[11pt]{amsart}
\usepackage{amsmath,amsthm}
\usepackage{graphics}
\usepackage{latexsym}
\usepackage{verbatim}
\usepackage{amsmath}
\usepackage{amsthm}
\usepackage{amssymb}
\usepackage{epsfig}
\usepackage{epstopdf}
\usepackage{subfig}
\usepackage{hyperref}
\usepackage[]{hyperref}
\usepackage{setspace}
\usepackage[noend]{algpseudocode}
\usepackage{algorithmicx,algorithm}
\hypersetup{urlcolor=blue, citecolor=red}
\usepackage{dsfont}
\usepackage[table]{xcolor}
\usepackage{multirow}
\usepackage{float}
\usepackage{tikz}
\usetikzlibrary{arrows, calc, shapes, positioning, decorations.pathreplacing}
\usepackage{comment}
\usepackage{MnSymbol}
\usepackage{diagbox}
\usepackage{multirow}

\ifpdf
  \DeclareGraphicsExtensions{.eps,.pdf,.png,.jpg}
\else
  \DeclareGraphicsExtensions{.eps}
\fi
\newtheorem{theorem}{Theorem}[section]
\newtheorem{proposition}[theorem]{Proposition}
\newtheorem{corollary}[theorem]{Corollary}
\newtheorem{lemma}[theorem]{Lemma}

\title[Block preconditioning methods for asymptotic preserving scheme]{Block preconditioning methods for asymptotic preserving scheme arising in anisotropic elliptic problems}

\author[Lingxiao Li and Chang Yang]{}


\email{lilingxiao@lsec.cc.ac.cn}
\email{yangchang@hit.edu.cn}

\begin{document}
\maketitle

\centerline{\scshape Lingxiao Li}
\medskip
{\footnotesize
    \centerline{Institute of Applied Physics and Computational Mathematics}
    \centerline{ Fenghaodong Road, Haidian District,
Beijing, 100094, China}
}

\medskip

\centerline{\scshape Chang Yang}
\medskip
{\footnotesize
    \centerline{School of Mathematics, Harbin Institute of Technology,}
    \centerline{No. 92 West Dazhi Street, Nangang District, 150001 Harbin, China}
}

\begin{abstract}
Efficient and robust iterative solvers for strong anisotropic elliptic equations are very challenging.
  In this paper a block preconditioning method is introduced to solve the linear algebraic systems of a
  class of micro-macro asymptotic-preserving (MMAP) scheme.
  MMAP method was developed by Degond {\it et al.} in 2012
  where the discrete matrix has a $2\times2$ block structure.
  By the approximate Schur complement a series of block preconditioners are constructed.
  We first analyze a natural approximate Schur complement that is the coefficient matrix of the original non-AP discretization.
  However it tends to be singular for very small anisotropic parameters.
  We then improve it by using more suitable approximation for boundary rows of the exact Schur complement.
  With these block preconditioners, preconditioned GMRES iterative method is developed to solve the discrete equations.  
  Several numerical tests show that block preconditioning methods can be a robust strategy with respect to
  grid refinement and the anisotropic strengths.
\end{abstract}


\tableofcontents

\section{Introduction}
\label{sec:introduction}

In this paper, the following anisotropic elliptic equation is considered:
\begin{equation}\label{eq:anisotropic_pb}
\left\{
\begin{array}{ll}
-\Delta_\bot \phi^\varepsilon - \frac{1}{\varepsilon}\Delta_{\|}\phi^{\varepsilon} = f^\varepsilon, & \text{in } \Omega,\, \\
\phi^\varepsilon = 0, & \text{on }\Gamma_D,\, \\
\nabla_\bot \phi^\varepsilon \cdot \mathbf{n} + \frac{1}{\varepsilon}
\nabla_{\|}\phi^\varepsilon \cdot \mathbf{n} = 0, &  \text{on }\Gamma_N,
\end{array}
\right.
\end{equation}
where $\Omega$ is an open set in $\mathbb{R}^2$ or  $\mathbb{R}^3$, its boundary is defined by $\partial \Omega= \bar{\Omega}\backslash\Omega$, $\mathbf{n}$ is exterior normal vector on the boundary  $\partial\Omega$. The parameter $\varepsilon$ is reciprocal of anisotropic strength, and it can be very small, {\it i.e.} $0<\varepsilon\ll 1$. Let $\mathbf{b}$ be a magnetic field, then we can divide the boundary $\partial\Omega$ into two parts, {\it i.e.} $\partial \Omega = \Gamma_D\cup\Gamma_N$, with
\begin{align*}
\Gamma_D & =\{x\in\partial\Omega \ |\ \mathbf{b}\cdot \mathbf{n} = 0\},\quad
\Gamma_N  =\{x\in\partial\Omega \ | \ \mathbf{b}\cdot \mathbf{n} \neq 0\}.
\end{align*}
Moreover, we define parallel or perpendicular differential operators as follows:
\begin{equation*}
\begin{array}{ll}
\nabla_{\|} \phi^\varepsilon = \mathbf{b}\otimes\mathbf{b}\nabla\phi^\varepsilon, & \nabla_{\bot} \phi^\varepsilon = (I - \mathbf{b}\otimes\mathbf{b})\nabla\phi^\varepsilon, \, \\
\Delta_\| \phi^\varepsilon = \nabla\cdot \nabla_{\|}\phi^\varepsilon, & \Delta_\bot \phi^\varepsilon = \nabla\cdot \nabla_{\bot}\phi^\varepsilon.
\end{array}
\end{equation*}
where $I$ is the identity matrix.
This anisotropic elliptic equation~\eqref{eq:anisotropic_pb} has important applications in the physics of magnetized plasma,
such as confinement fusion~\cite{Meiss2003Plasma}, plasma of ionosphere~\cite{Besse2004A}, plasma propulsion~\cite{Meier2011}, {\it  etc}.

The problem~\eqref{eq:anisotropic_pb} is well-posed for any $\varepsilon>0$.
Indeed let $\Omega$ be a convex polygon and $f^\varepsilon\in L^2(\Omega)$,
then we have a unique solution $\phi^\varepsilon\in H^2_{D}(\Omega)$  of the problem~\eqref{eq:anisotropic_pb}, with
\begin{equation*}
H^2_{D}(\Omega) = \{v\in H^2(\Omega):\, v|_{\Gamma_D} = 0\}.
\end{equation*}

Formally, if $\varepsilon$ goes to 0, we get a degenerate problem:
\begin{equation}\label{eq:degenerate_pb}
\left\{
\begin{array}{ll}
-\Delta_{\|}\phi^{0} = 0, & \text{in } \Omega,\, \\
\phi^0 = 0, & \text{on }\Gamma_D,\, \\
\nabla_{\|}\phi^0 \cdot \mathbf{n} = 0, &  \text{on }\Gamma_N.
\end{array}
\right.
\end{equation}
The problem~\eqref{eq:degenerate_pb} is ill-posed due to loss of uniqueness of solution, since any function belongs to set 
\[H^2_{D,\|} = \{v\in H^2(\Omega):\, v|_{\Gamma_D} = 0,\, \nabla_\|v=0\text{ in }\Omega\}\]
is solution of~\eqref{eq:degenerate_pb}.
In other words, the solution of~\eqref{eq:degenerate_pb} is in kernel of differential operator $\nabla_\|$.

We then define a limit problem
\begin{equation}\label{eq:limit_pb}
\left\{
\begin{array}{ll}
-\Delta_{\bot}\bar\phi = f^0, & \text{in } \Omega,\, \\
\bar\phi = 0, & \text{on }\Gamma_D,\, \\
\nabla_{\bot}\bar\phi \cdot \mathbf{n} = 0, &  \text{on }\Gamma_N,\, \\
\nabla_{\|}\bar\phi = 0, &  \text{in } \Omega.
\end{array}
\right.
\end{equation}
This limit problem~\eqref{eq:limit_pb} has unique solution $\bar\phi \in H^2_{D,\|}$ thanks to the constraint $\nabla_{\|}\bar\phi = 0  \text{ in } \Omega$.

Solving directly the anisotropic elliptic equation with classic numerical methods may involve severe numerical pollution~\cite{yang2020Preserving}. One possible way to alleviate such pollution of anisotropic elliptic or parabolic problems is to use high order methods, such as~\cite{2008APSDPPBP6040K} for fourth-order finite volume method and~\cite{2020Giorgiani} for high order finite element method. Moreover, Chac\'on et al. in~\cite{L2014An} proposed an asymptotic-preserving semi-Lagrangian algorithm for time-dependent anisotropic heat conduction equations. 

In the last decade, another class of  asymptotic preserving  (AP) schemes were developed for the problem~\eqref{eq:anisotropic_pb}.
Those AP schemes, on the one hand, are equivalent to the problem~\eqref{eq:anisotropic_pb} when $\varepsilon>0$; on the other hand, they coincide with the limit problem~\eqref{eq:limit_pb} as $\varepsilon\to0$. Therefore, the AP schemes are well-defined for all values of $\varepsilon$. The first AP scheme was given for aligned case ($\mathbf{b}=(0,1)^T$) in~\cite{degond_asymptotic_2009}. Some extensions were made for aligned case, namely Duality Based methods~\cite{Christophe2013Efficient,2018Iterative}. Other methods focused on general non-aligned case, such as Micro-Macro AP (MMAP) scheme~\cite{DEGOND20122724}, field line integration method~\cite{tang_asymptotic_2017, wang_uniformly_2018}, two field iterated method~\cite{deluzet_two_2019, yang_numerical_2019}, AP scheme based on first-order system method~\cite{2021Li}. In most of these methods, the key point is to introduce an auxiliary variable, then we have an AP system. To solve the AP system numerically,
one must solve a linear algebraic system which has more degree of freedoms than original systems, where direct solvers
may not be practical in large scale simulation. Thus the main purpose of this paper is to develop efficient and robust iterative methods.

In this paper we propose to use GMRES method~\cite{Saad1986} to solve the linear algebraic systems. To improve the robustness and efficiency of
the iterative solver, the key ingredient is efficient preconditioning. The coefficient matrix generated from the MMAP scheme used in this paper
has a typical $2\times 2$ block structure. Motivated by the block preconditioning methods from Navier-Stokes equations in~\cite{nsref1}, where
the key point is approximate Schur complement, we attempt to devise efficient approximate Schur complement of the MMAP system in this paper.
We firstly give a natural approximation from the continuous operator level. However this natural approximate Schur complement may be
ill conditioned if $\varepsilon$ is very small. We then replace its boundary rows by the ones of exact Schur complement. It is surprising
that the GMRES method converges only in 2 iterations in our numerical examples.
Therefore we find the core point is how to approximate the boundary rows of exact Schur complement.
Based on this fact, we present some improvements in Section 3.2.

This paper is organized as follows: In Section 2, the Micro-Macro Asymptotic Preserving (MMAP) schemes introduced in~\cite{DEGOND20122724} is briefly revised. In Section 3, several block preconditioning methods for the MMAP schemes are described. In Section 4, a number of numerical comparison of those  block preconditioning methods are given, including the aligned case and non-aligned case. The preconditioned methods for the MMAP scheme without inflow boundary condition is also discussed here. Finally, we finish this paper with some conclusion and perspectives.

\section{Review of the Micro-Macro Asymptotic Preserving schemes}
\label{sec:MMAPscheme}
In this part, we briefly revise the  Micro-Macro Asymptotic Preserving (MMAP) schemes introduced in~\cite{DEGOND20122724}. For clarity, let us omit the subscript $\varepsilon$. Moreover, let us focus on an anisotropic elliptic equation defined in a unit square, {\it e.g.} $\bar{\Omega}=[0,1]\times[0,1]$. Moreover, let us denote $\Omega= (0,1)\times(0,1)$, $\Gamma_D = \{0,1\}\times[0,1]$ and $\Gamma_N = [0,1]\times\{0,1\}$.
Thus the anisotropic elliptic equation is given as follows
\begin{equation}\label{eq:anisotropic_pb2}
\left\{
\begin{array}{ll}
-\Delta_\bot\phi - \frac{1}{\varepsilon}\Delta_{\|} \phi = f,& \text{in }\Omega,\\[3mm]
\phi=0,&\text{on }\Gamma_D,\\[3mm]
\mathbf{n}\cdot \nabla_\bot \phi + \frac{1}{\varepsilon}\mathbf{n}                                                                                                                                                                                                                                                                                                                                                       \cdot \nabla_\| \phi=0,&\text{on } \Gamma_N.
\end{array}
\right.
\end{equation}
The key point of the MMAP scheme in~\cite{DEGOND20122724} is to introduce an auxiliary function $q$, with the following ansatz
\begin{equation}\label{eq:auxiliary}
\nabla_\|\phi = \varepsilon\nabla_\| q.
\end{equation}
Taking another time of parallel derivative on~\eqref{eq:auxiliary}, we have
\begin{equation}\label{eq:auxiliary2}
\Delta_\|\phi = \varepsilon\Delta\| q, \text{  in }\Omega.
\end{equation}
Injecting~\eqref{eq:auxiliary2} into the first equation of~\eqref{eq:anisotropic_pb2}, we get
\begin{equation*}
-\Delta_\bot\phi - \Delta_\| q = f, \text{  in }\Omega.
\end{equation*}
Similarly, putting the ansatz~\eqref{eq:auxiliary} into the boundary conditions in~\eqref{eq:anisotropic_pb2}, we can obtain the corresponding boundary conditions. Finally, the MMAP scheme is written as
\begin{equation}\label{eq:MMAP}
\left\{
\begin{array}{ll}
-\Delta_\bot\phi - \Delta_\| q = f,& \text{  in }\Omega,\\[3mm]
\Delta_\|\phi - \varepsilon\Delta_\| q = 0,& \text{  in }\Omega,\\[3mm]
\phi=0,&\text{  on }\Gamma_D,\\[3mm]
\mathbf{n}\cdot \nabla_\bot \phi + n\cdot \nabla_\| q=0,&\text{  on } \Gamma_N,\\[3mm]
\mathbf{n}\cdot \nabla_\|\phi - \varepsilon n\cdot \nabla_\| q = 0,& \text{  on }\Gamma_N.
\end{array}
\right.
\end{equation}

Now, let us give two remarks on the MMAP scheme~\eqref{eq:MMAP}. The first one is that the MMAP scheme has a unique solution $\phi$ for any $\varepsilon\geq0$ provided that the parallel derivative of $\phi$ is vanishing on the ends of the magnetic field lines. Indeed, provided that $\varepsilon>0$, the  MMAP scheme~\eqref{eq:MMAP} is equivalent to the anisotropic elliptic equation~\eqref{eq:anisotropic_pb2}, which is well-posed for $\varepsilon>0$ (see~\cite{DEGOND20122724}), thus we have unique solution $\phi$. Formally letting $\varepsilon\to0$ in~\eqref{eq:MMAP}, we have
\begin{equation}\label{eq:deltaparaphivanish}
\Delta_\|\phi = 0\text{  in }\Omega.
\end{equation}
Integrating~\eqref{eq:deltaparaphivanish} along magnetic field line and thanks to $\nabla_\|\phi$ vanishing at one end of magnetic field line, we find that
\begin{equation*}
\nabla_\|\phi = 0\text{  in }\bar\Omega,
\end{equation*}
which means that $\phi$ is constant along  magnetic field line, thus $\phi\equiv\bar{\phi}$, where we denote $\bar{\phi}$ as the average of $\phi$ along a magnetic field line. At last, integrating the first equation of~\eqref{eq:MMAP} along a magnetic field line, we obtain a unique solution $\bar{\phi}$, which is also the unique solution for the MMAP scheme~\eqref{eq:MMAP}. Notice that the constraint the parallel derivative of $\phi$  vanishing on  the ends of the magnetic field lines is equivalent to say that the  magnetic field lines cross the boundary $\Gamma_N$.

The second remark is that the MMAP scheme~\eqref{eq:MMAP} does not determine a unique $q$. Indeed, if $q$ satisfies~\eqref{eq:MMAP}, then $q+q_0$ also satisfies~\eqref{eq:MMAP} for any $q_0$ in the kernel of the operator $\nabla_\|$. To determine a unique $q$, in~\cite{DEGOND20122724}, the inflow boundary condition was added on one end of magnetic field lines (denoted as  $\Gamma_{\text{in}}\subset\Gamma_{N}$ where $\mathbf{b}\cdot\mathbf{n}<0$), that is $q=0$ on $\Gamma_{\text{in}}$. Moreover, notice that along a magnetic field line, one boundary condition of $\mathbf{n}\cdot \nabla_\|\phi - \varepsilon \mathbf{n}\cdot \nabla_\| q = 0$ on $\Gamma_N$ is tedious, since one can easily verify this fact by integrating the second equation of~\eqref{eq:MMAP} along a magnetic field line. As a result, on one end of magnetic field lines, the boundary condition $\mathbf{n}\cdot \nabla_\|\phi - \varepsilon \mathbf{n}\cdot \nabla_\| q = 0$ can be substituted by the inflow boundary condition $q=0$.

\section{Block preconditioning methods}\label{sec:blockprecon}
In this part, we will introduce block preconditioning methods for discrete MMAP scheme of Section~\ref{sec:MMAPscheme}. Let us reformulate the MMAP scheme as follows
\begin{equation}\label{eq:MM1}
\left\{
\begin{array}{ll}
 - \Delta_\| q-\Delta_\bot\phi = f,& (x,z)\in(0,1)\times(0,1),\\[3mm]
q=0,&(x,z)\in\{0,1\}\times[0,1],\\[3mm]
q=0,&(x,z)\in [0,1]\times\{0\}, \\[3mm]
\mathbf{n}\cdot \nabla_\| q + \mathbf{n}\cdot \nabla_\bot \phi =0,&(x,z)\in [0,1]\times\{1\},
\end{array}
\right.
\end{equation}
and
\begin{equation}\label{eq:MM2}
\left\{
\begin{array}{ll}
\varepsilon\Delta_\| q - \Delta_\|\phi  = 0,& (x,z)\in(0,1)\times(0,1),\\[3mm]
\phi=0,&(x,z)\in\{0,1\}\times[0,1],\\[3mm]
\mathbf{n}\cdot \nabla_\| q + \mathbf{n}\cdot \nabla_\bot \phi =0,&(x,z)\in [0,1]\times\{0\},\\[3mm]
\varepsilon \mathbf{n}\cdot \nabla_\| q - \mathbf{n}\cdot \nabla_\|\phi  = 0,& (x,z)\in[0,1]\times\{1\}.
\end{array}
\right.
\end{equation}
In the sequel, two cases will be considered:
\begin{itemize}
\item Aligned case: the magnetic field line aligned with $z$-axis, {\it i.e.} $\mathbf{b}=(0,1)^T$.

\item Non-aligned case: the magnetic field line does not aligned with Cartesian coordinates.
\end{itemize}
In Appendix~\ref{sec:FDD},  the finite difference methods for the MMAP scheme~\eqref{eq:MM1}-\eqref{eq:MM2} are given for both aligned case and non-aligned case. In any case, the linear system takes a matrix form as
\begin{equation}\label{eq:MMAPlinsyst}
\begin{pmatrix}
A_3 & B \\
A_2 & A_1
\end{pmatrix}
\begin{pmatrix}
Q \\
\Phi
\end{pmatrix}
=
\begin{pmatrix}
F \\
0
\end{pmatrix}.
\end{equation}
This block matrix can be factorized as
\begin{equation*}
\mathbf{A} =
\begin{pmatrix}
A_3 & B \\
A_2 & A_1
\end{pmatrix}
=
\begin{pmatrix}
A_3 & 0 \\
A_2 & A_1-A_2 A_3^{-1}B
\end{pmatrix}
\begin{pmatrix}
I & A_3^{-1}B \\
0 & I
\end{pmatrix}.
\end{equation*}
The so called Schur complement is
\begin{equation}\label{eq:schurcomplement}
E = A_1-A_2 A_3^{-1}B.
\end{equation}
Denote the lower block triangular matrix as follows
\begin{equation*}
\mathbf{P}
=
\begin{pmatrix}
A_3 & 0 \\
A_2 & {E}
\end{pmatrix}.
\end{equation*}
Using $\mathbf{P}$
as preconditioner, typical Krylov subspace methods such as  GMRES~\cite{Saad1986} will converge in two iterations~\cite{Murphy-ref}~\cite{2001Efficient}
for any grid size and any value of $\varepsilon$.
However, the price to pay is to invert the matrix $A_3$, which is too costly when discretization mesh size $N$ is large enough.

Now the central question is how to find a suitable approximation $S$ of the exact Schur complement $E = A_1 - A_2 A_3^{-1} B$ with a cheap cost.
More precisely, let us consider the following block preconditioner
\begin{equation*}
\mathbf{P}_S
=
\begin{pmatrix}
A_3 & 0 \\
A_2 & S
\end{pmatrix},
\end{equation*}
its inverse is
\begin{align*}
\mathbf{P}_S^{-1}
=
\begin{pmatrix}
A_3^{-1} & 0 \\
-S^{-1}A_2A_3^{-1} & S^{-1}
\end{pmatrix}.
\end{align*}
 Now applying the preconditioner to $\mathbf{A}$, we get
\begin{equation*}
\mathbf{P}_S^{-1}\mathbf{A}=
\begin{pmatrix}
I & A_3^{-1}B \\
0 &  S^{-1}E
\end{pmatrix},
\end{equation*}
which means the eigenvalues of $\mathbf{P}_S^{-1}\mathbf{A}$ are 1 and $\lambda$, where $\lambda$ is   eigenvalue of $ S^{-1}E$. If we could construct a suitable approximation for the exact Schur complement $E$, then the eigenvalues of $\mathbf{P}_S^{-1}A$ will be clustered tightly. As a result one can expect that the Krylov iterative methods have a rapid convergence.
\subsection{Approximate Schur complement}\label{sec:ApproximateSchurcomplement}

In this subsection we give the basic ideas to develop the approximate Schur complement for the MMAP system matrix $\mathbf{A}$ in~\eqref{eq:MMAPlinsyst}. We express the coefficient matrix of~\eqref{eq:MMAPlinsyst} in continuous operator form
\begin{equation*}
\begin{pmatrix}
A_3 & B \\
A_2 & A_1
\end{pmatrix}
\approx
\begin{pmatrix}
-\Delta_\| & -\Delta_\bot \\
\varepsilon\Delta_\| & -\Delta_\|
\end{pmatrix}.
\end{equation*}
Formally, one can verify that the underlying continuous operator for the exact Schur complement $E = A_1 - A_2 A_3^{-1} B$ is
\begin{equation}\label{eq:continuousschur}
-\Delta_\| - \varepsilon\Delta_\|(-\Delta_\|)^{-1}(\Delta_\bot) = -\Delta_\| - \varepsilon\Delta_\bot,
\end{equation}
which is precisely the original PDE operator for variable $\phi$ in~\eqref{eq:anisotropic_pb2}. Thus a natural approximation matrix for $E$ is the discrete coefficient matrix for~\eqref{eq:continuousschur} with zero boundary condition on $\Gamma_D$ and zero Neumann boundary condition on $\Gamma_N$. We denote such \textbf{natural approximate Schur complement} as $S_1$. However when $\varepsilon\to0$, due to zero parallel derivative $\mathbf{n}\cdot\nabla_\|\phi$, one can deduce that $S_1$ tends to be singular.
The performance of such block preconditioner with $S_1$ will deteriorate.
For this assertion, one can see numerical results of aligned case in Table~\ref{tab:gmresiternum_aligned}.

Now we want to give a new approximation from the algebraic level. Let us denote the diagonal part of $A_3$ by $D$,
which may offer some useful spectral information of $A_3$. Hence, substituting  $A_3$ by $D$, we obtain a new approximation for $E$
\begin{equation}
S_2 = A_1 - A_2D^{-1}B.
\end{equation}
The matrix $S_2$ contains the spectral information of $A_3$ in the algebraic level and will show extra advantage when $\varepsilon$ becomes very small as indicated by our numerical results in section~\ref{sec:num}.
We remark that the proposed two approximate Schur complement $S_1$ and $S_2$ can be easily constructed in any space dimension and be conveniently accessed in parallel computing environment.
\vspace{3mm}

\textbf{Remark 1.} The natural approximate Schur complement $S_1$ needs more explanation here.
Firstly it performs rather well for large $\varepsilon$.
Moreover, in our numerical experiments we find that the iteration numbers
of GMRES using $S_1$ is also stable in non-aligned case. But this is not in aligned case. In non-aligned case, we numerically observe 
that the condition number of $S_1$ scales as follows
\[\textsf{cond}(S_1)\sim \frac{1}{\varepsilon + \beta^2 h^2}\frac{1}{h^2}\]
where $\beta$ is a parameter of magnetic field in \eqref{eq:magneticfield1}.
Therefore if $\beta$ is not very small and the grid size $h$ is not very fine,
this "bad approximation" of non-AP discretization can be
viewed as a good preconditioning for the MMAP discretization.

\subsection{A further modification}

Above we have constructed two approximate Schur complement. One is motivated from the continuous PDE level and the other mainly from the algebraic level. In this subsection we combine the two ideas to give further modification for the approximation. 

Recall that the operator for $A_1$ is $-\Delta_\|$  and for $B$ is $-\Delta_\bot$, so we can also use matrix $A_1+\varepsilon B$ to approximate operator~\eqref{eq:continuousschur}. Unfortunately, on the one hand, the matrix $A_1$ may have  zero rows on $[0,1]\times\{ 0\}$, since $\mathbf{n}\cdot\mathbf{b}^{\bot}$ ($\mathbf{b}^\bot$ is orthogonal to $\mathbf{b}$) may vanish on  $[0,1]\times\{ 0\}$, and on the other hand, matrix $B$ has no boundary conditions on $[0,1]\times\{ 0\}$. This fact means that  $A_1+\varepsilon B$ may be singular due to null rows on bottom boundary $[0,1]\times\{ 0\}$. To remedy this, we replace such null rows in $A_1+\varepsilon B$  by the corresponding rows in $S_2$. Then we obtain a third approximate Schur complement $S_3$.


In order to assess the influence of the approximation of the boundary rows in $S_3$,
we construct an extra expense approximation $S_4$ using exact inverse of $A_3$.
In detail we replace null boundary rows of $A_1+\varepsilon B$  on $[0,1]\times\{ 0\}$
using the corresponding rows of the exact Schur complement $A_1-A_2A_3^{-1}B$
instead of the ones of $S_2$. Here $S_4$ is just used for performance comparison
and we do not recommend this use in practice.

We also propose two additional approximate Schur complements. The first one, denoted as $S_5$,
is a combination of $S_1$ and $S_3$, that is  the rows corresponding to  $[0,1]\times\{ 0\}$ of $S_5$ are simply sum of the rows of $S_1$ and $S_3$. From this, we expect to incorporate the both advantages of  $S_1$ and $S_3$.
The second one is only valid for aligned case with the discretization given in~\ref{sec:discrete_aligned},
where we expect to construct a symmetric positive definite approximate Schur complement.
For this, we impose a Robin boundary condition as follows
\begin{equation*}
\phi-h\partial_z\phi=0,\quad (x,z)\in(0,1)\times\{0\}.
\end{equation*}
The corresponding discretization is
\begin{equation*}
2\phi_{i,1} - \phi_{i,2}=0, \quad i\in [2,n_x-1].
\end{equation*}
Let us denote this approximate Schur complement by $S_6$.

Now let us give a short summarize. We first construct three approximate Schur complements  $S_\alpha,\alpha=1,2,3$ from different point of view.
We also use an expense $S_4$ for performance comparison. Finally, we give two additional ones $S_5$ and $S_6$ for improved performance.
\vspace{3mm}

\textbf{Remark 2.} In the appendix, we prove that for aligned case the exact Schur complement $E$ is precisely $S_4$. Our numerical tests
also seem to indicate this because the GMRES solver converges in only 2 steps for both aligned case and non-aligned case
and for any grid size and anisotropic parameters. This fact tells us that if we can find good approximation of boundary rows of
exact Schur complement then we can find good approximate Schur complement.

\section{Numerical results}\label{sec:num}
This section is devoted to show performance of the proposed block preconditioner of Section~\ref{sec:blockprecon}.
We notice that to apply $\mathbf{P}_{S}, $ as a preconditioner, in each step of GMRES iteration,
a following solve of  triangular block linear system is needed
\begin{equation*}
\begin{pmatrix}
A_3 & 0 \\
A_2 & S_{\alpha}
\end{pmatrix}
\begin{pmatrix}
e_1\\ e_2
\end{pmatrix}
=
\begin{pmatrix}
r_1\\ r_2
\end{pmatrix},
\end{equation*}
which is equivalent to solve two linear systems successively as
\begin{eqnarray}
 A_3 e_1&=&r_1,\label{eq:pcsystem1}\\
 S_{\alpha}e_2&=&r_2-A_2e_1.\label{eq:pcsystem2}
\end{eqnarray}
\subsection{Setup}
\begin{figure}

\begin{center}
\begin{tabular}{cc}
\includegraphics[width=5cm]{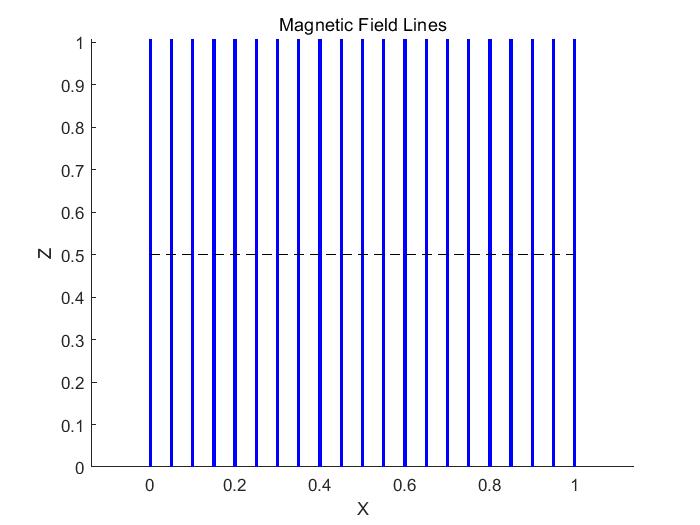} & \includegraphics[width=5cm]{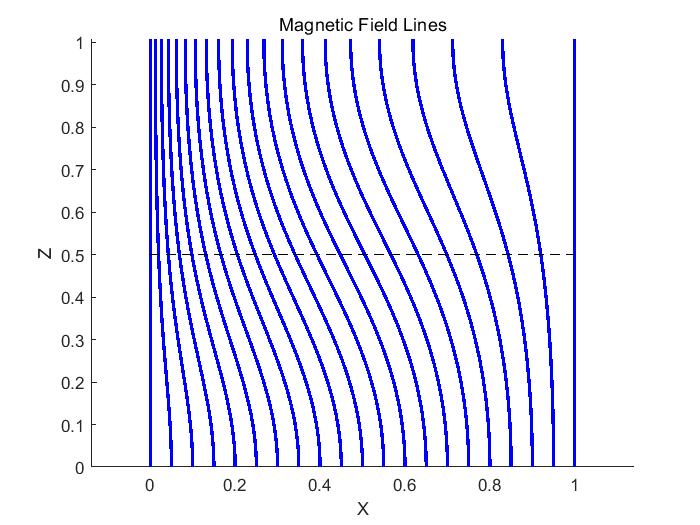} \\
(a) $\beta=0$ & (b) $\beta=2$
\end{tabular}
\end{center}
\caption{\label{fig:magneticfield}Magnetic field lines defined in~\eqref{eq:magneticfield1}.}
\end{figure}

In this section, we consider a fixed magnetic field~\cite{DEGOND20122724,deluzet_two_2019}
\begin{equation}\label{eq:magneticfield1}
B = \begin{pmatrix}
\pi\beta(x^2-x)\sin(\pi z)\\
\beta(2x-1)\cos(\pi z)+\pi
\end{pmatrix}.
\end{equation}
Thus we have
\begin{equation}\label{eq:magneticfield2}
b = \frac{B}{|B|}.
\end{equation}
The manufactured solution $\phi$ is given by
\begin{equation}
\phi(x,z) = \sin(m(\pi x + \beta (x^2 - x)\cos(\pi z))) + \varepsilon \cos(2\pi z)\sin(\pi x).
\end{equation}
Below, following parameters are used, $m=4$, $\beta=0$ or $2$, $\varepsilon\in[10^{-20},1]$. The magnetic field lines defined in~\eqref{eq:magneticfield1} are illustrated in Figure~\ref{fig:magneticfield}, where $\beta=0$ corresponds to aligned case and 
$\beta=2$ for non-aligned case.
In the sequel, without special explanation, the generalized minimal residual algorithm (GMRES)~\cite{Saad1986} method with no restart will be used as iterative solver for solving the linear system~\eqref{eq:MMAPlinsyst}. The relative tolerance is set to be 1.0e-6, and the maximum iterative number is 100. Moreover, all the resolution of the MMAP scheme have reach desired precision, thus we will not show errors analysis in this section.

\subsection{Aligned case}\label{sec:numres_aligned}

In this part, we consider aligned case, where discretization is described in Appendix~\ref{sec:discrete_aligned}.

\begin{figure}

\begin{center}
\begin{tabular}{ccc}
\includegraphics[width=3.8cm]{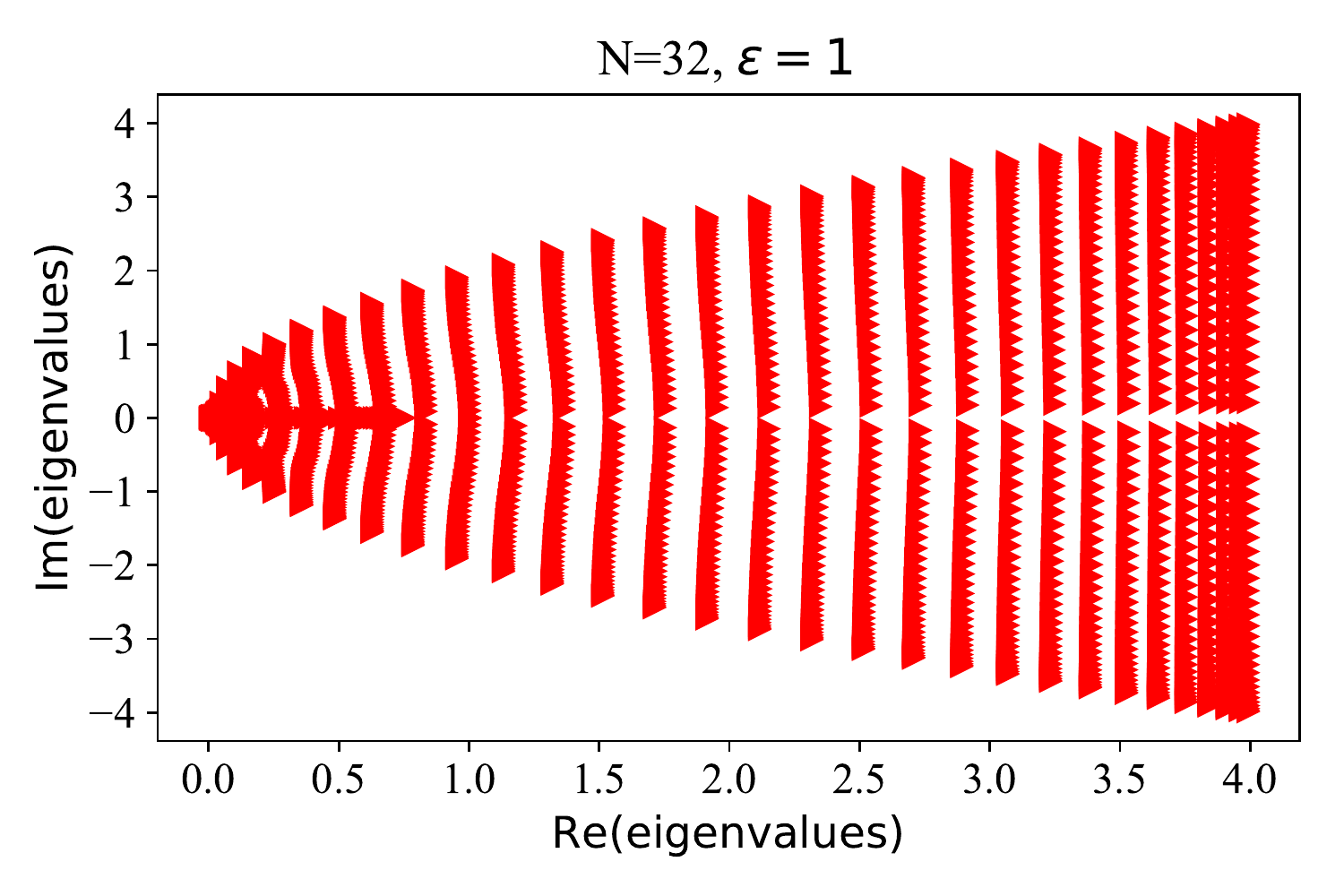} &
 \includegraphics[width=3.8cm]{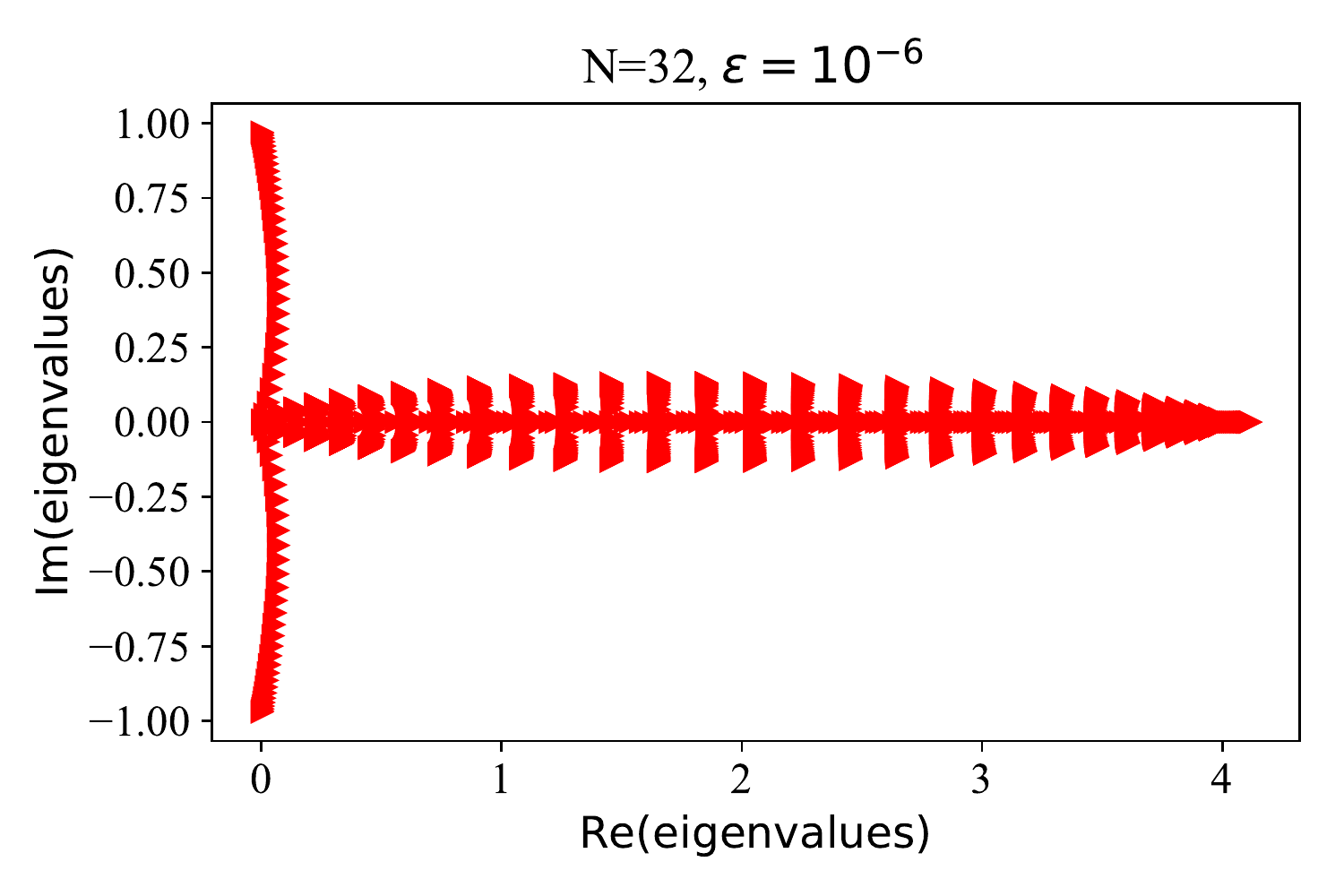} &
\includegraphics[width=3.8cm]{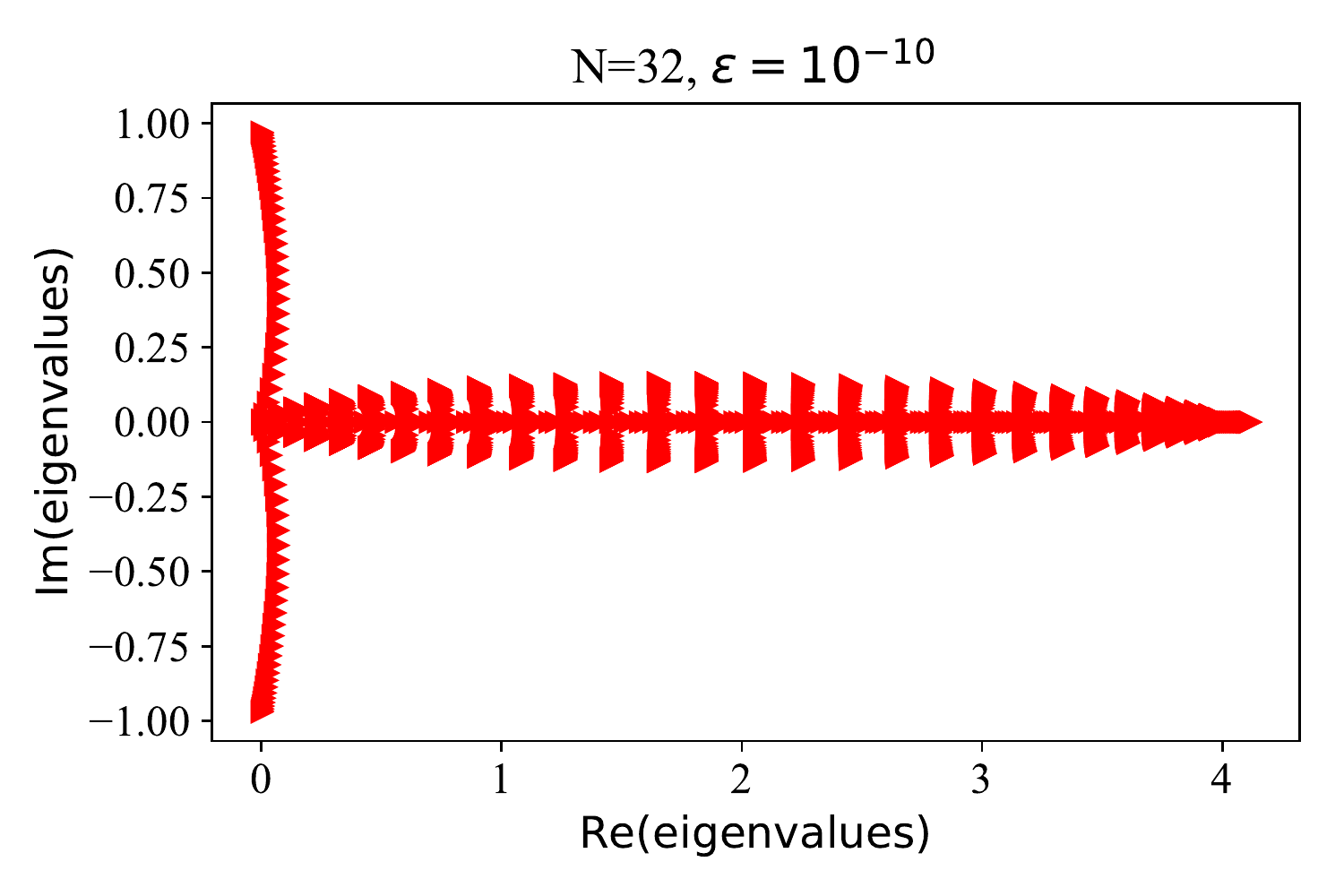} \\
\includegraphics[width=3.8cm]{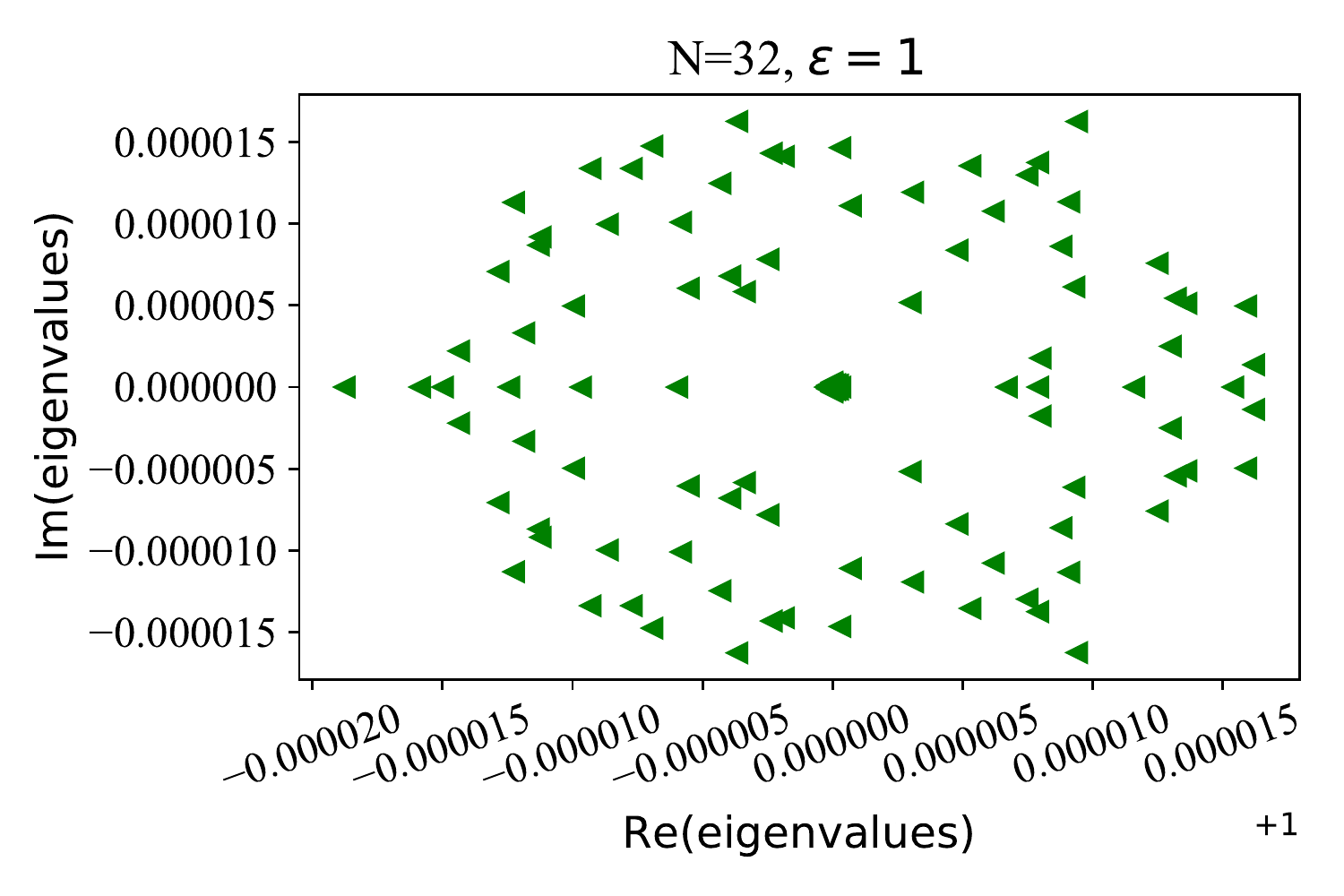} & \includegraphics[width=3.8cm]{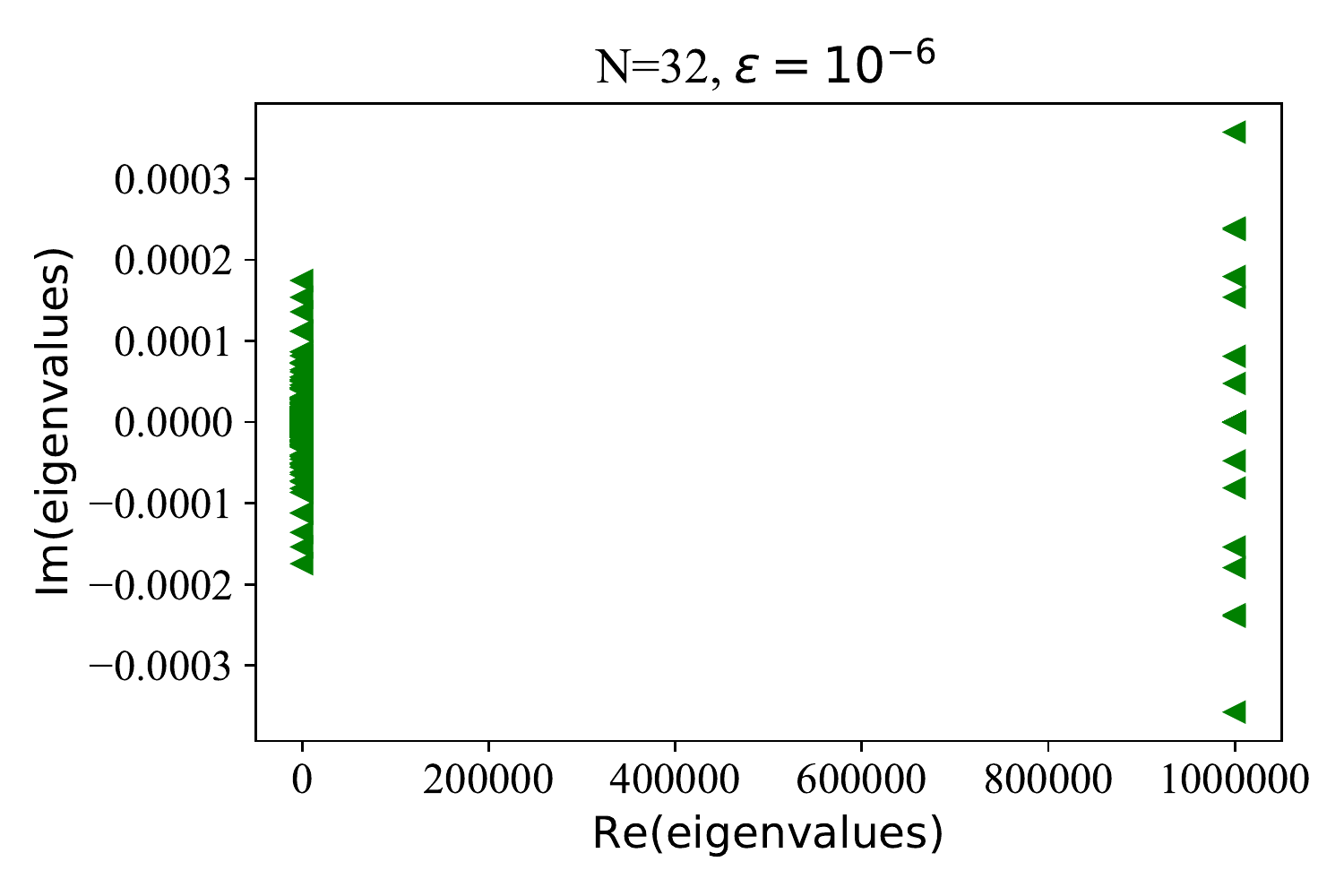} &
\includegraphics[width=3.8cm]{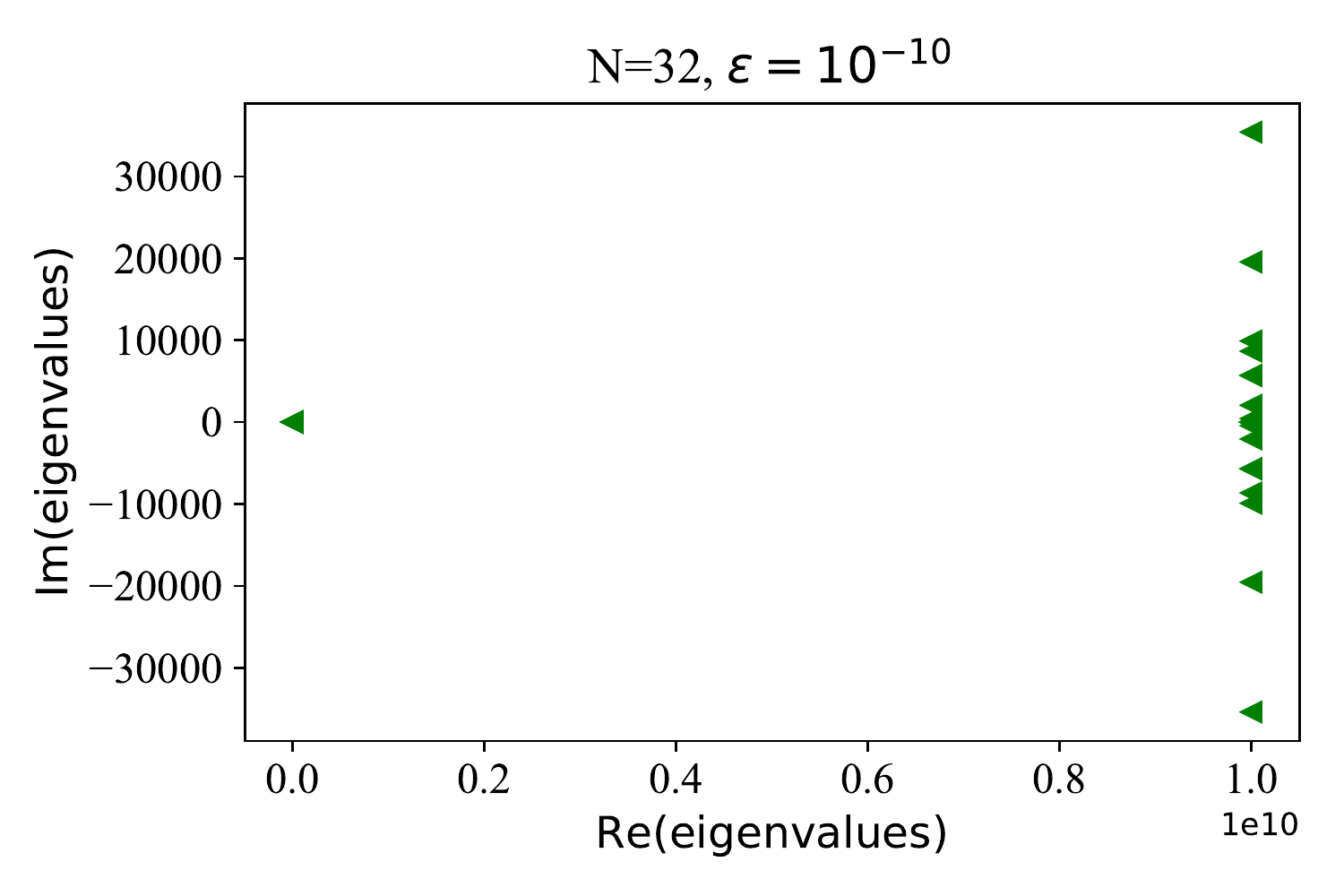} \\
\includegraphics[width=3.8cm]{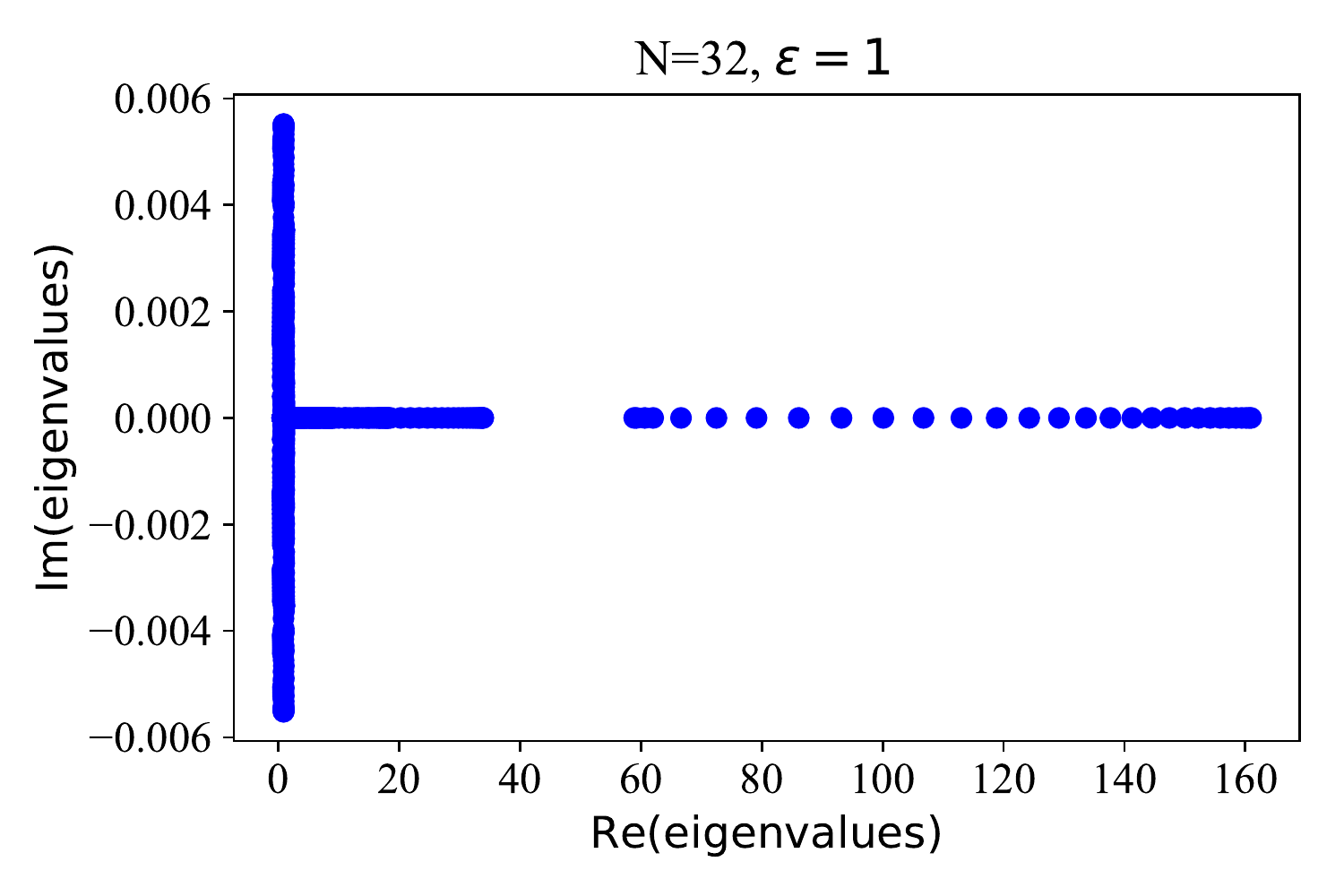} & \includegraphics[width=3.8cm]{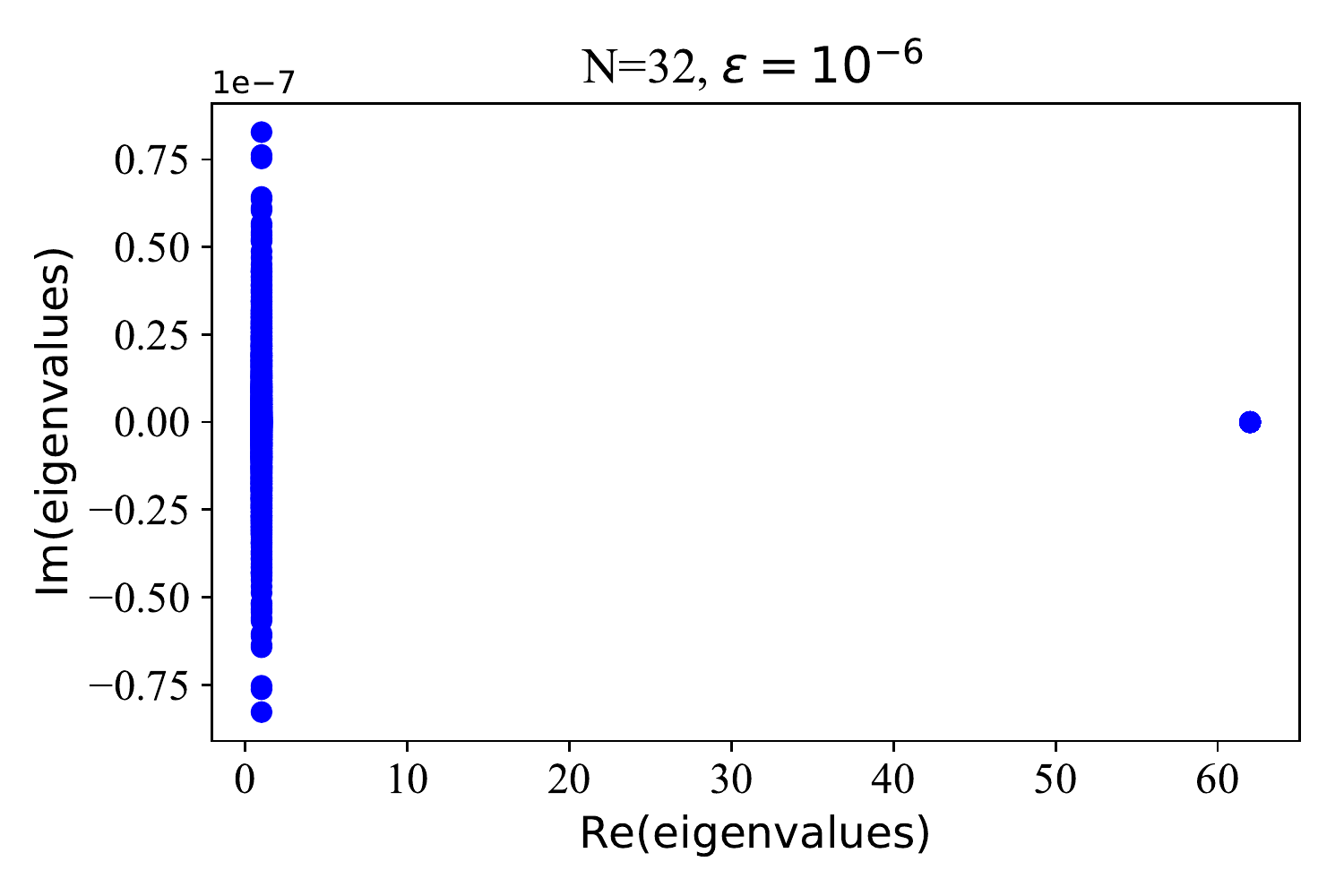} &
\includegraphics[width=3.8cm]{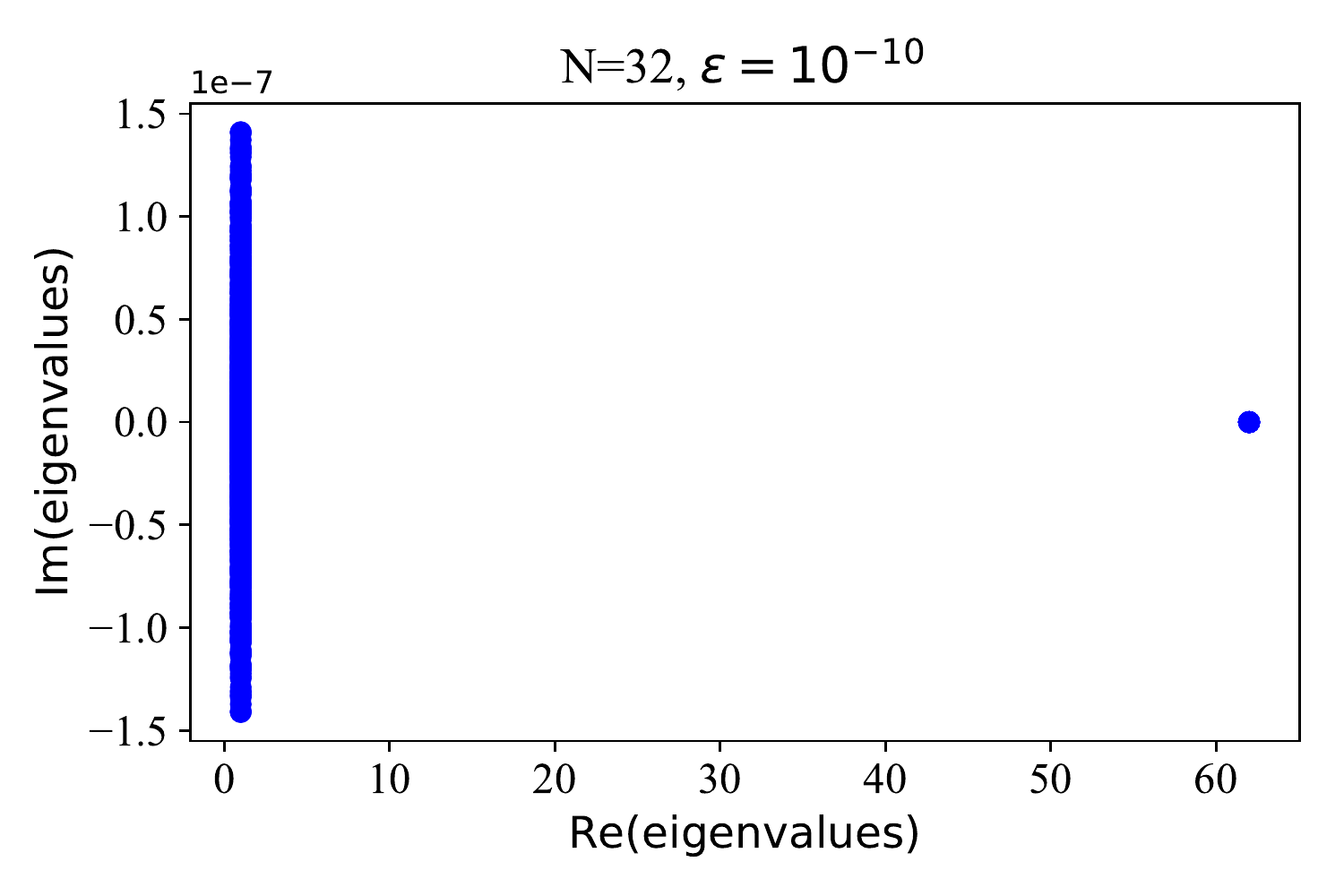} \\
\includegraphics[width=3.8cm]{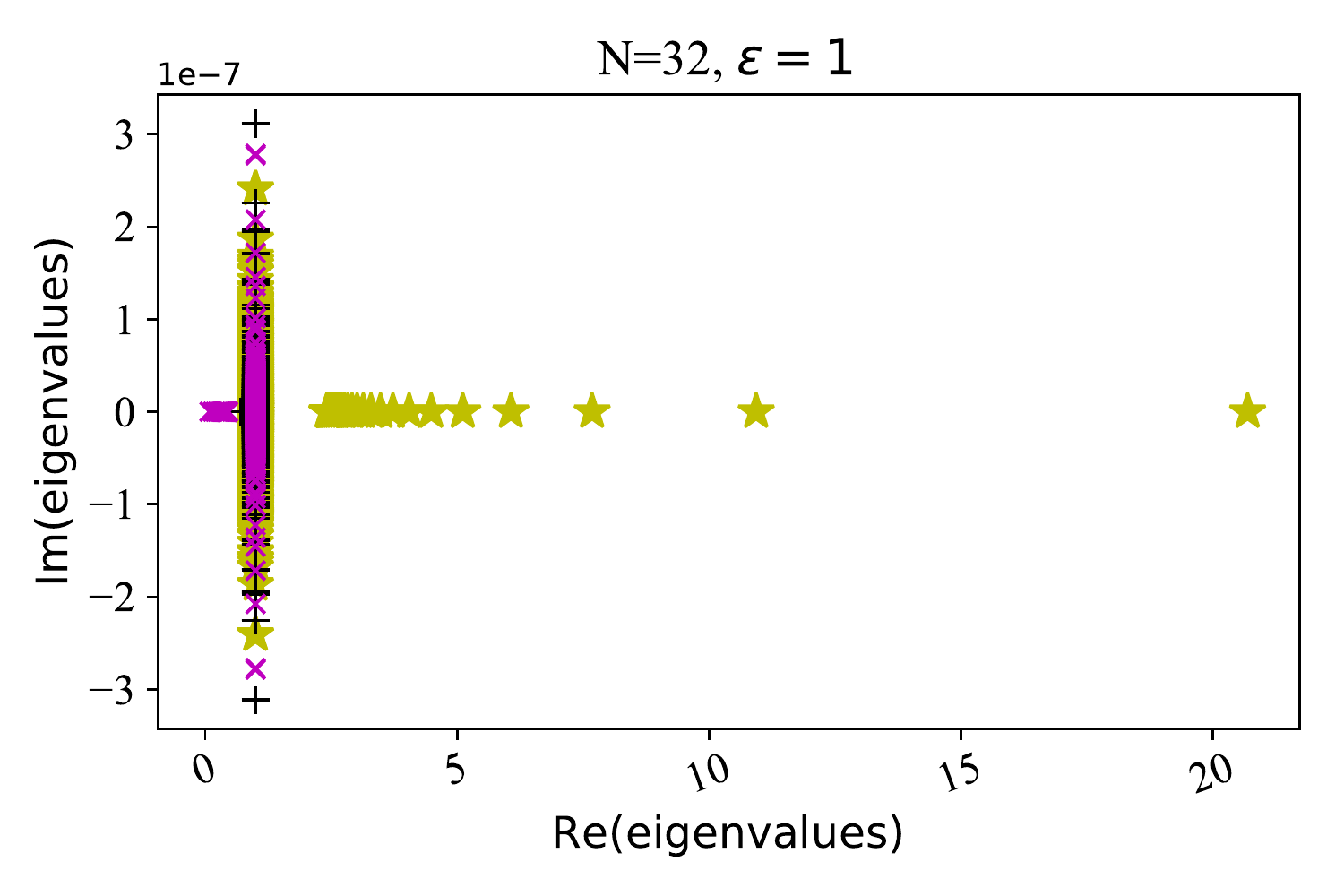} & \includegraphics[width=3.8cm]{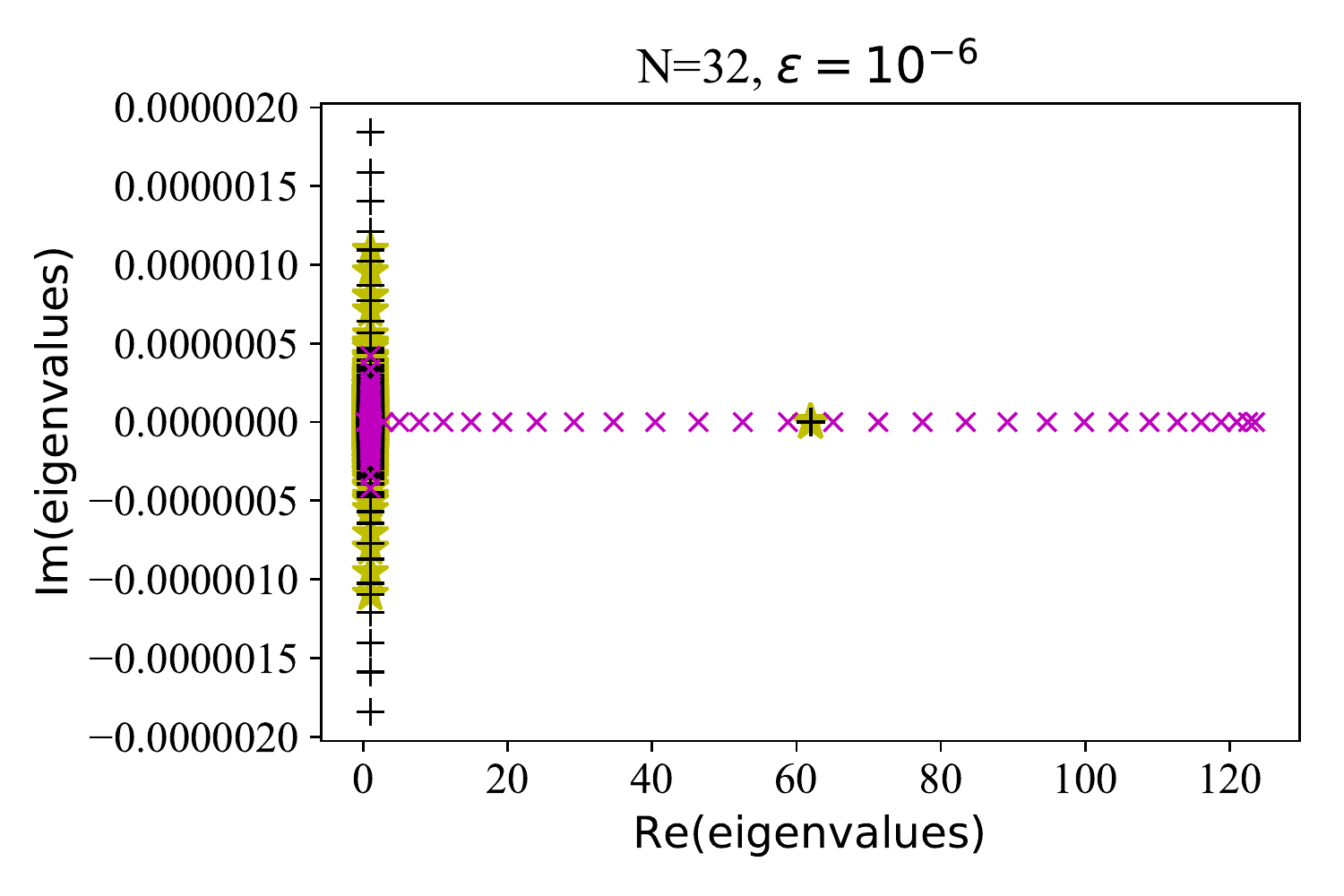} &
\includegraphics[width=3.8cm]{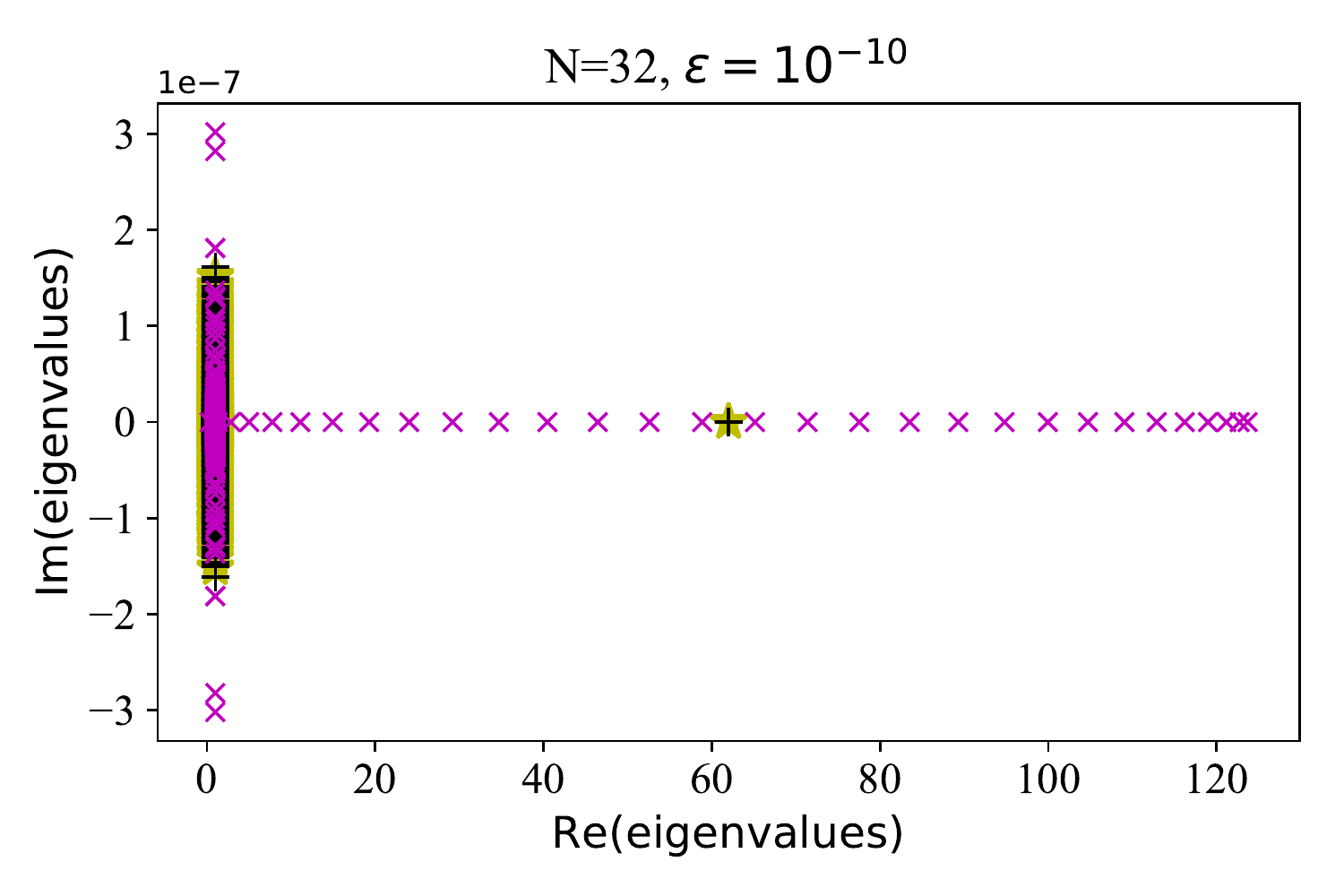}
\end{tabular}
\end{center}
\caption{\label{fig:aligned_eigenvalues}
  Aligned case. Eigenvalues for the original system $\mathbf{A}$ and the preconditioned system $\mathbf{P}_S^{-1}\mathbf{A}$.
  Top row: unpreconditioned eigenvalues (red $\blacktriangleright$'s);
  second row: $S_1$ (green $\blacktriangleleft$'s);
  third row: $S_2$ (blue $\bullet$'s);
  bottom row: $S_3, S_5, S_6$ (yellow $\star$'s, black $+$'s, magenta $\times$'s respectively).
  $N$ stands for mesh size.
}
\end{figure}

Let us first show graphically the eigenvalues of different matrices,  which is reported in Figure~\ref{fig:aligned_eigenvalues}. We observe that the eigenvalues of the unpreconditioned matrix $\mathbf{A}$ have both large distribution in real axis and image axis, and this similar phenomena is valid for all values of $\varepsilon$. On the contrary, for preconditioned matrix $\mathbf{P}^{-1}_S\mathbf{A}$, we expect eigenvalues are tightly close to positive real axis and are as clustered as possible. Ideally, we hope these eigenvalues tend to 1. For the case with $S_1$, we see the extent of eigenvalues is  proportional to $\frac{1}{\varepsilon}$. This is due to the fact that there is no cross derivatives in aligned case, thus the discretizations of $\Delta_\|$ and  $\Delta_\bot$ can be made separately, {\it i.e.} $\Delta_\|$  is discretized aligned with $z-$axis and the one of $\Delta_\bot$ is aligned with $x-$axis.  The eigenvalues of $S_1$ are real and the maximum range can be estimated as~\cite{2018Iterative}
\begin{equation*}
\max(\textsf{eig}(S_1)) \sim  \frac{1}{\varepsilon h^2}.
\end{equation*}
For the case with $S_2$, we see the eigenvalues already becomes more clustered than the case with $S_1$, especially for small $\varepsilon$, {\it e.g.} $\varepsilon<10^{-6}$.
However, for $\varepsilon=1$, the eigenvalues are not tightly clustered to the positive real axis, and this phenomena can be observed more significantly on refined mesh.
This is because when $\varepsilon$ is big, the solutions of original systems are close to solutions of 
MMAP scheme. In such case $S_1$ ia a very good approximation of $E$.
However this is not the case with very small $\varepsilon$.

For the case with $S_\alpha,\alpha=3,5,6$, we see that the eigenvalues are tightly close to the positive real axis for all values of $\varepsilon$. More precisely, we find that for $\varepsilon=1$, the eigenvalues for cases with $S_5$ and $S_6$ are very similar, and they are better than the case with $S_3$. However, for small $\varepsilon$, the eigenvalues for cases with $S_3$ and $S_5$ become more similar, and they are better than the case with $S_6$. As a result, the case with $S_5$ seems to be the best in the comparison shown in Figure~\ref{fig:aligned_eigenvalues}.
Notice that in Figure~\ref{fig:aligned_eigenvalues}, the eigenvalues of the case with $S_4$ are not shown, it is because $S_4$ is
nothing but the exact Schur complement $E$ (see the detail in Appendix~\ref{sec:S4}), thus the eigenvalues of this case are all equal to 1.


Now let us solve the linear system~\eqref{eq:MMAPlinsyst} with GMRES solver.
At each iteration we need a resolution of the
linear systems~\eqref{eq:pcsystem1}-\eqref{eq:pcsystem2}. This resolution is divided into two steps: the first step is to factorize the matrices $A_3$ and $S_\alpha$, then the second step is to solve~\eqref{eq:pcsystem1}-\eqref{eq:pcsystem2} thanks to the factorized matrices obtained in the first step. Due to the finite difference discretization in Appendix~\ref{sec:discrete_aligned}, we know that the matrices $A_3$ and $S_6$ are symmetric positive definite (SPD), we thus factorize them by   the Cholesky decomposition. For the other $S_\alpha,\alpha=1,\dots,5$, we just apply the LU decomposition.


The iteration numbers of the GMRES method are reported in Table~\ref{tab:gmresiternum_aligned}, where we compare among different mesh sizes and anisotropy strength. Overall, the results we got are consistent with the previous analysis of eigenvalues. For the case with $S_1$,
at most 3 iterations are needed for convergence of the GMRES method when $\varepsilon\geq10^{-10}$. However when $\varepsilon$ is vanishing, the GMRES method does not converge anymore. For the case with  $S_2$, the preconditioner works well for small $\varepsilon$, but the convergence rate is quite slow for large $\varepsilon$. For the case with   $S_3$, it combines the advantages of cases with $S_1$ and $S_2$.
We observe that the iterations do not increase when refining mesh and the iterations decrease for vanishing $\varepsilon$. For the case with $S_4$, just two iterations are needed since $S_4$ is exact Schur complement. For the case with $S_5$, we obtain exact table as for the case with $S_3$. Finally, for the case with $S_6$, the results are very similar as the one for the case with  $S_3$ or  $S_5$, except for $\varepsilon=10^{-6}$ and $N=256$, where the GMRES solver converges in 6 iterations.


\begin{table}[!h]
\caption{\label{tab:gmresiternum_aligned} Iteration numbers of GMRES method with approximate Schur complements $S_\alpha$, $\alpha=1,\dots,6$.}
\begin{center}
(a) Case with  $S_1$

\begin{tabular}{ccccccc}
\hline
$\varepsilon$ & 1 & $10^{-1}$ & $10^{-2}$ & $10^{-6}$ & $10^{-10}$ & $10^{-20}$\\
\hline
$N = 32$ & 3 & 3 & 3 & 3 & 1 &  $\times$  \\
$N = 64$ & 3 & 3 & 3 & 3 & 2 &  $\times$  \\
$N = 128$ & 3 & 3 & 3 & 3 & 2 &  $\times$  \\
$N = 256$ & 3 & 3 & 3 & 3 & 2 & $\times$ \\
\hline
\end{tabular}
\vspace{3mm}

(b) Case with  $S_2$

\begin{tabular}{ccccccc}
\hline
$\varepsilon$ & 1 & $10^{-1}$ & $10^{-2}$ & $10^{-6}$ & $10^{-10}$ & $10^{-20}$\\
\hline
$N = 32$ & 14 & 6 & 5 & 2 & 2 & 2 \\
$N = 64$ & 18 & 7 & 5 & 2 & 2 & 2 \\
$N = 128$ & 33 & 11 & 4 & 2 & 2 & 2 \\
$N = 256$ & 67 & 11 & 4 & 2 & 2 & 2 \\
\hline
\end{tabular}
\vspace{3mm}

(c) Case with  $S_3$

\begin{tabular}{ccccccc}
\hline
$\varepsilon$ & 1 & $10^{-1}$ & $10^{-2}$ & $10^{-6}$ & $10^{-10}$ & $10^{-20}$\\
\hline
$N = 32$ & 4 & 4 & 4 & 3 & 2 & 2 \\
$N = 64$ & 4 & 4 & 4 & 3 & 2 & 2 \\
$N = 128$ & 4 & 4 & 4 & 3 & 2 & 2 \\
$N = 256$ & 4 & 4 & 4 & 4 & 2 & 2 \\
\hline
\end{tabular}
\vspace{3mm}

(d) Case with  $S_4$

\begin{tabular}{ccccccc}
\hline
$\varepsilon$ & 1 & $10^{-1}$ & $10^{-2}$ & $10^{-6}$ & $10^{-10}$ & $10^{-20}$\\
\hline
$N = 32$ & 2 & 2 & 2 & 2 & 2 & 2 \\
$N = 64$ & 2 & 2 & 2 & 2 & 2 & 2 \\
$N = 128$ & 2 & 2 & 2 & 2 & 2 & 2 \\
$N = 256$ & 2 & 2 & 2 & 2 & 2 & 2 \\
\hline
\end{tabular}
\vspace{3mm}

(e) Case with  $S_5$

\begin{tabular}{ccccccc}
\hline
$\varepsilon$ & 1 & $10^{-1}$ & $10^{-2}$ & $10^{-6}$ & $10^{-10}$ & $10^{-20}$\\
\hline
$N = 32$ & 4 & 4 & 4 & 3 & 2 & 2 \\
$N = 64$ & 4 & 4 & 4 & 3 & 2 & 2 \\
$N = 128$ & 4 & 4 & 4 & 3 & 2 & 2 \\
$N = 256$ & 4 & 4 & 4 & 4 & 2 & 2 \\
\hline
\end{tabular}
\vspace{3mm}

(f) Case with  $S_6$

\begin{tabular}{ccccccc}
\hline
$\varepsilon$ & 1 & $10^{-1}$ & $10^{-2}$ & $10^{-6}$ & $10^{-10}$ & $10^{-20}$\\
\hline
$N = 32$ & 4 & 4 & 4 & 3 & 2 & 2 \\
$N = 64$ & 4 & 4 & 4 & 3 & 2 & 2 \\
$N = 128$ & 4 & 4 & 4 & 2 & 2 & 2 \\
$N = 256$ & 4 & 4 & 4 & 6 & 2 & 2 \\
\hline
\end{tabular}
\vspace{3mm}

\end{center}
\end{table}

\begin{figure}

\begin{center}
\begin{tabular}{cc}
\includegraphics[width=5cm]{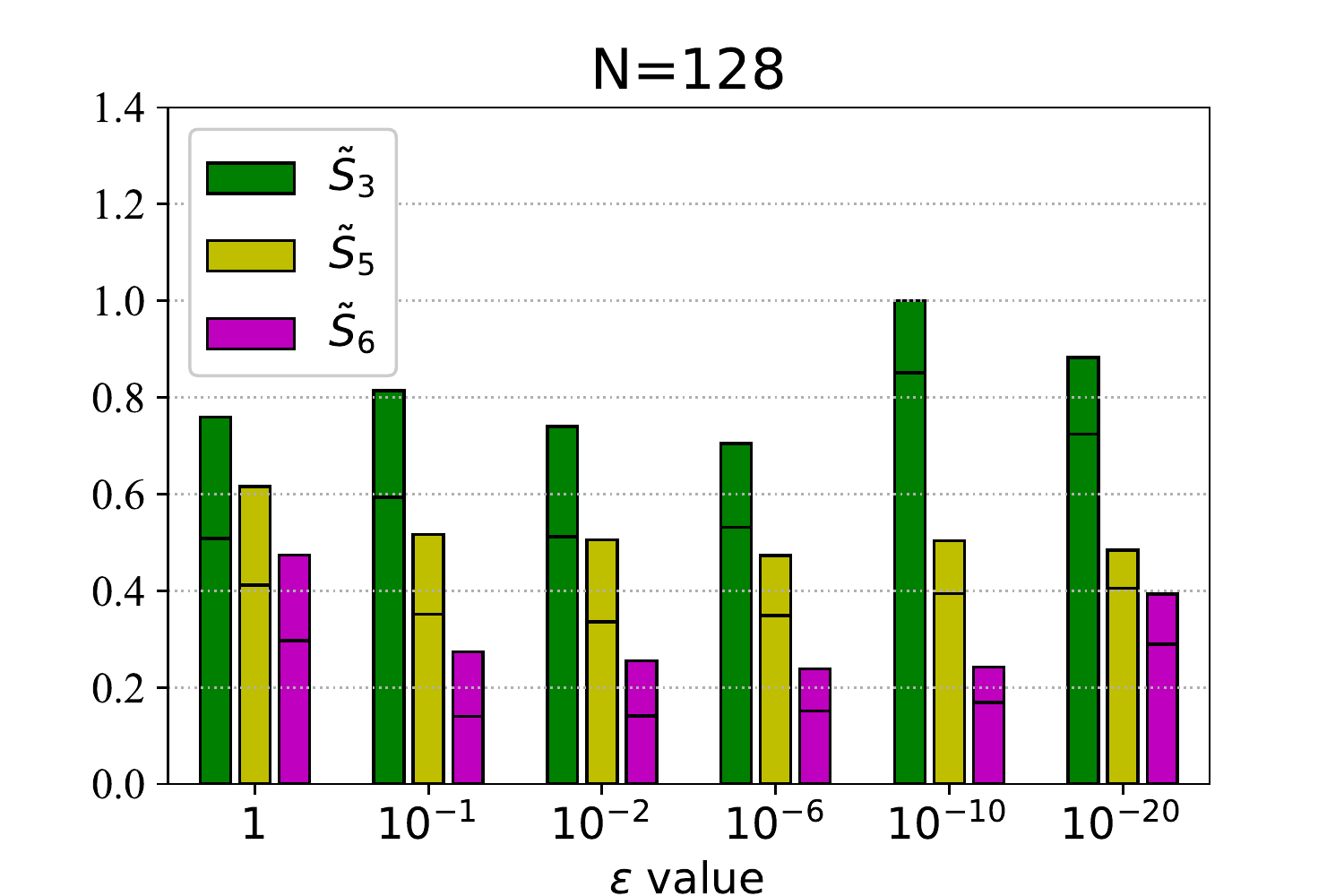} & \includegraphics[width=5cm]{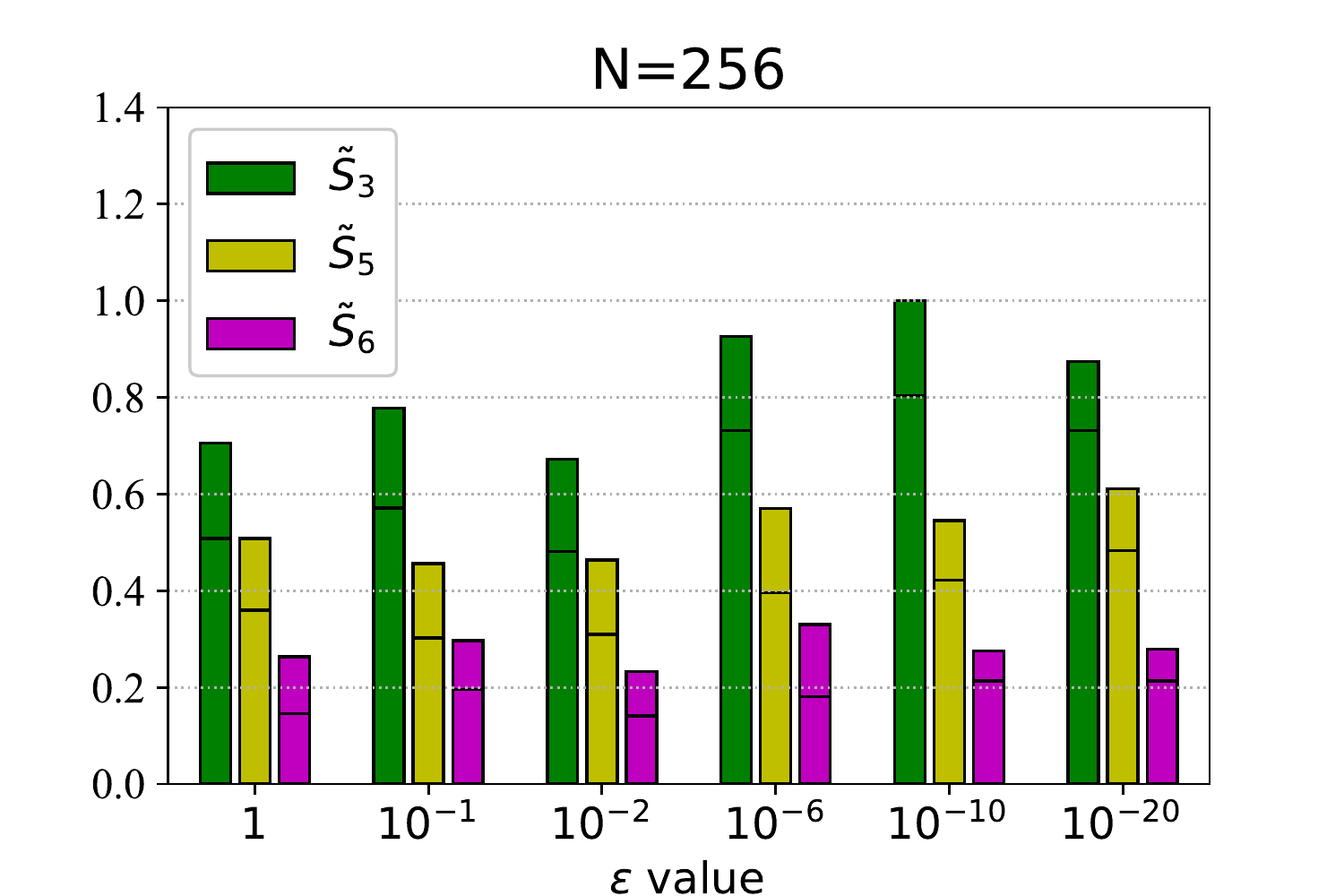}
\end{tabular}
\end{center}
\caption{\label{fig:aligned_directtime}Aligned case. The computational time for solving preconditioned linear system. The presented time is rescaled by the maximum value from each group of comparison. $\tilde{S}_\alpha,\alpha=3,5,6$ stands for different approximate Schur complement. Each bar is decomposed into two portions: bottom one is  time  of matrix factorizations and top one is for resolution time. $N$ stands for mesh size.}
\end{figure}

The computational time for solving preconditioned linear system is reported in Figure~\ref{fig:aligned_directtime}. Firstly, we see the  matrix factorization takes most of computational time, especially for very fine mesh, factorization time is much more than that for GMRES resolution. Secondly, Cholesky decomposition takes less than half of computational time of LU decomposition. This agrees with well-known results, that is complexity of Cholesky decomposition is $\mathcal{O}\left(\frac{1}{3}n^3\right)$ ($n$ is the size of matrix), while the one of LU decomposition is  $\mathcal{O}\left(\frac{2}{3}n^3\right)$.
Therefore, we conclude that  the resolution for the case with approximate Schur complement $S_6$ is more efficient than the cases with the other approximate Schur complements  $S_\alpha,\alpha=1,2,3,5$ for the aligned case.
Note that finite difference discretization given in Appendix~\ref{sec:discrete_nonaligned} can also deal with aligned case just by setting $\beta=0$. We find that for the case with $S_\alpha,\alpha=1,\dots,5$, we almost obtain the same results as presented in this subsection. However, since in this case $S_6$ is no longer symmetric, we can not benefit the advantages of SPD.


\subsection{Non-aligned case}\label{sec:numres_nonaligned}
\begin{figure}

\begin{center}
\begin{tabular}{ccc}
\includegraphics[width=3.8cm]{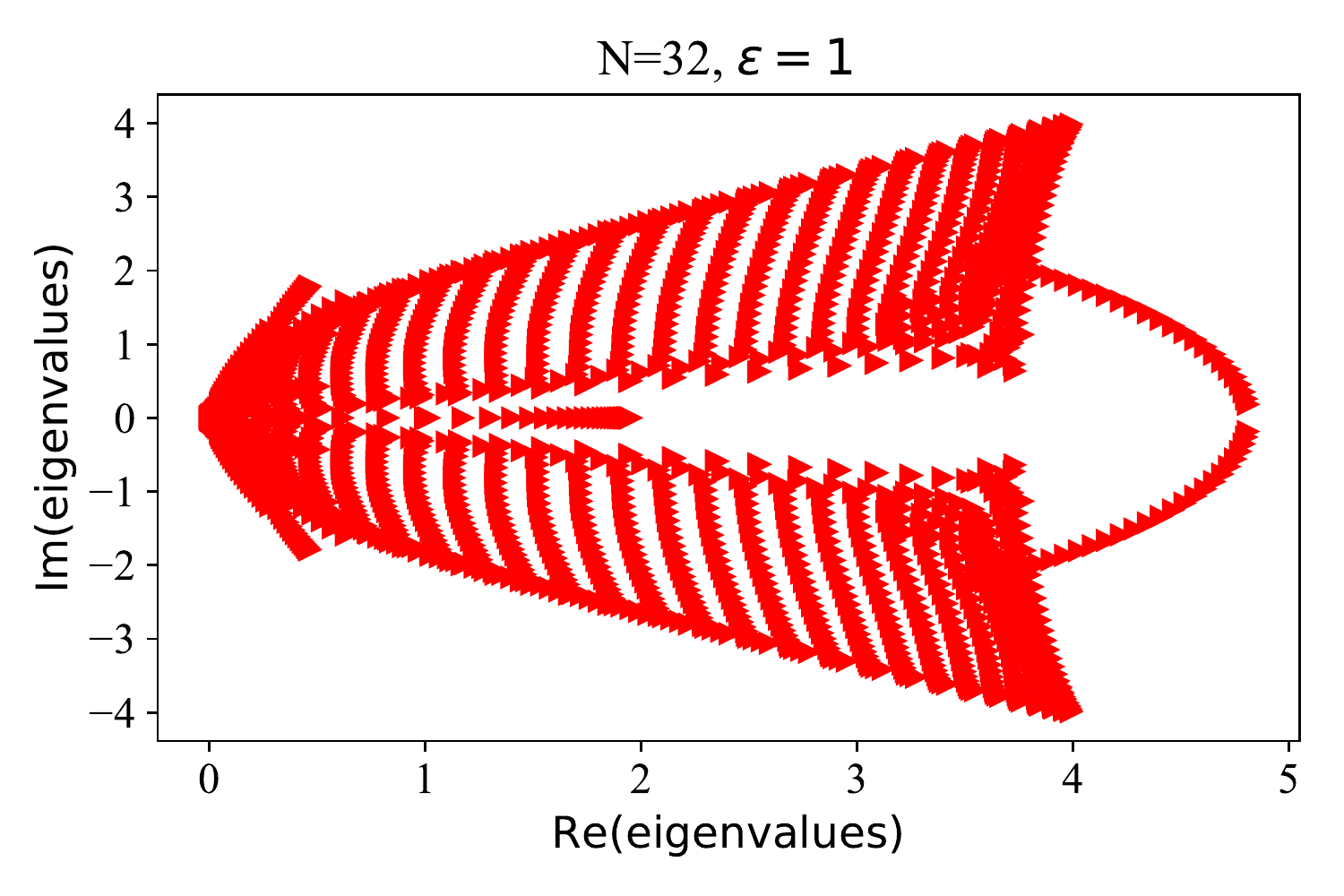} & \includegraphics[width=3.8cm]{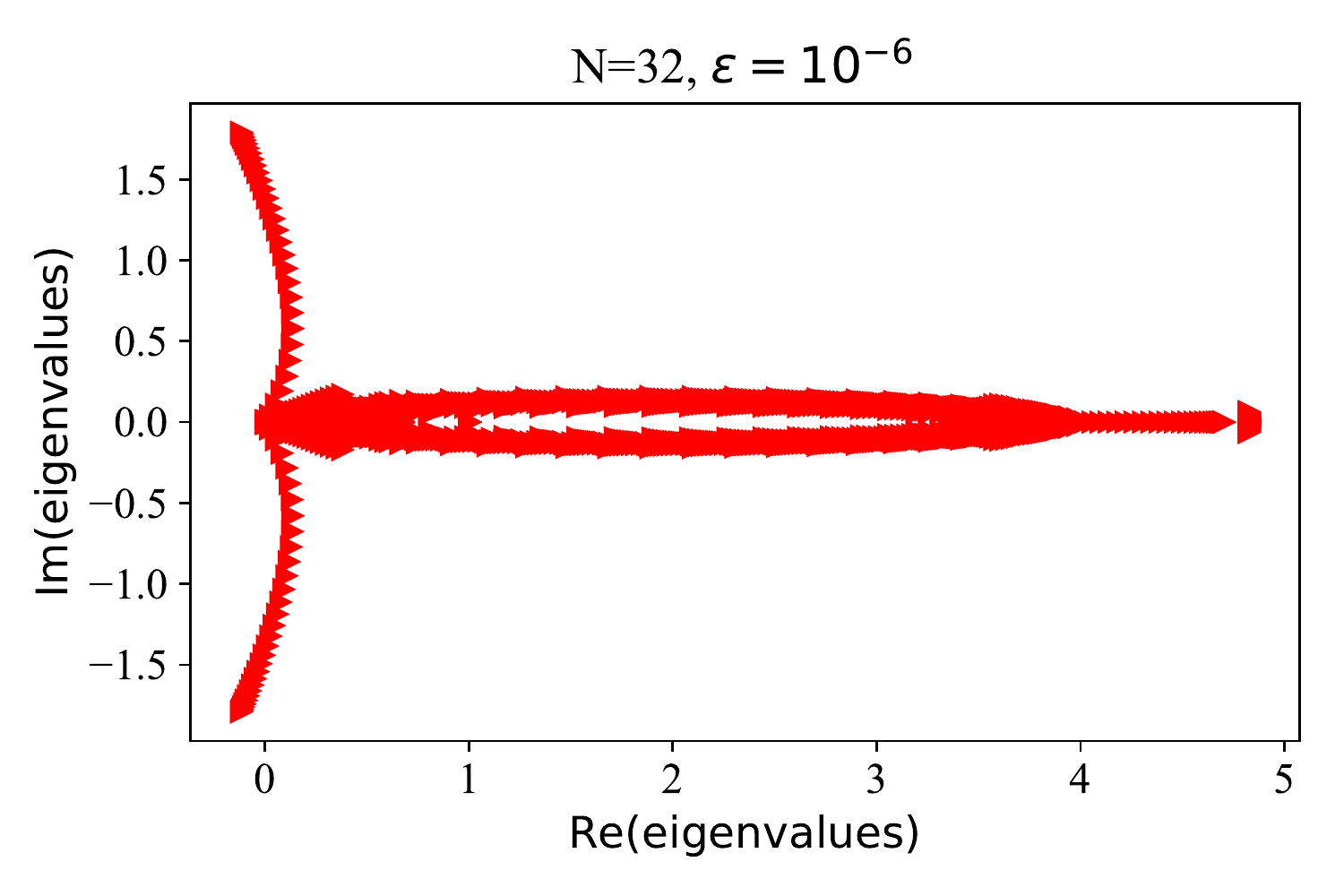} &
\includegraphics[width=3.8cm]{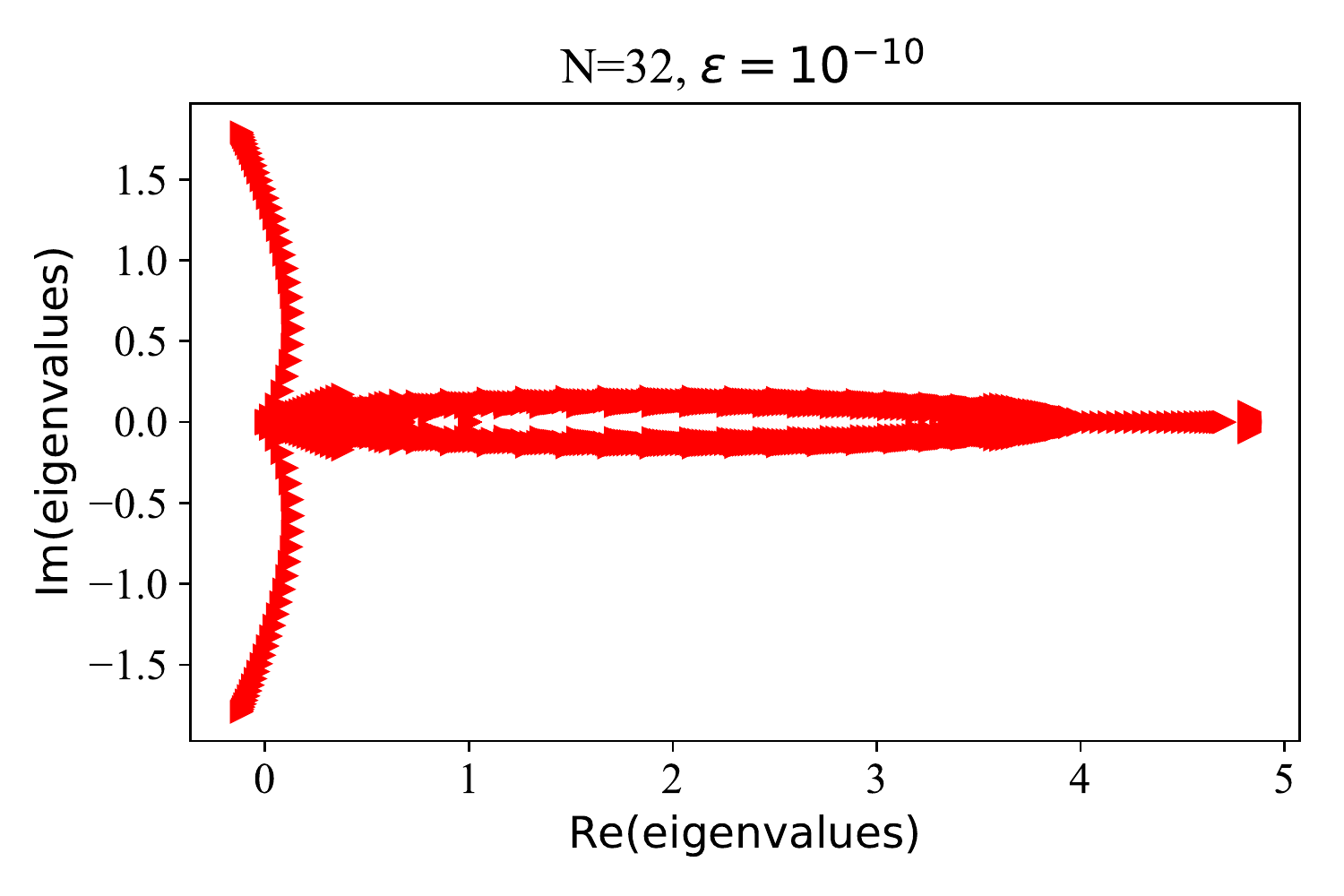} \\
\includegraphics[width=3.8cm]{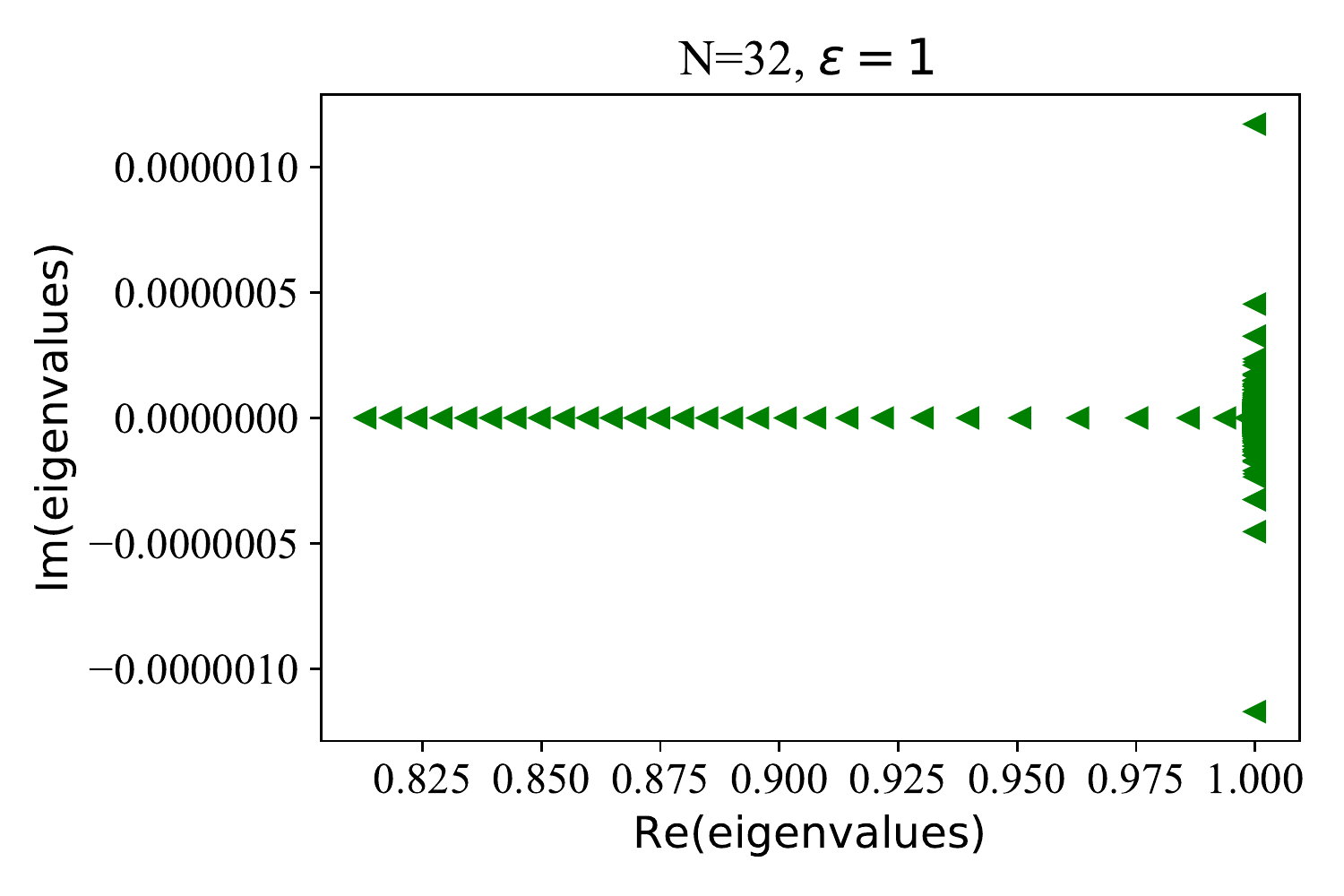} & \includegraphics[width=3.8cm]{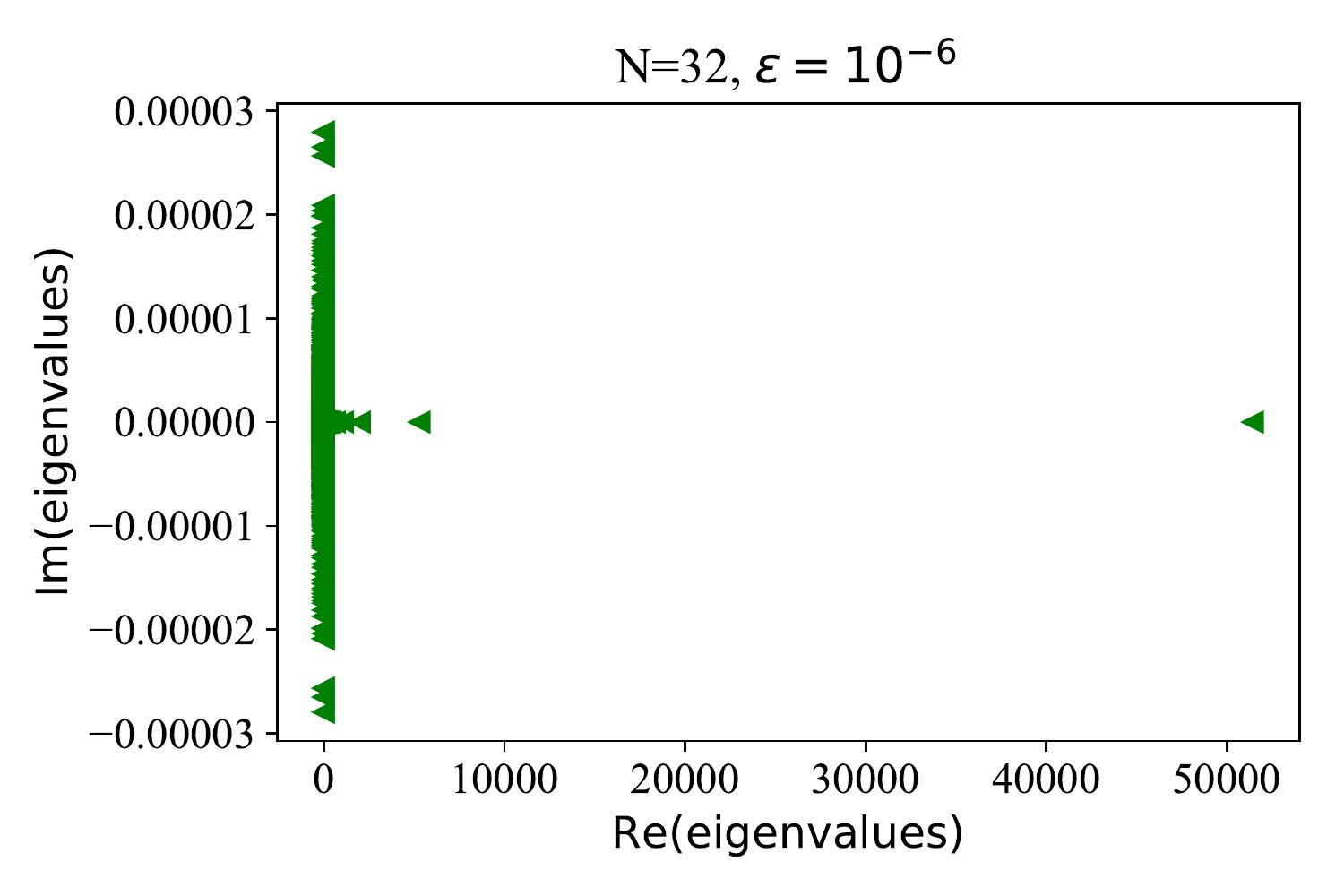} &
\includegraphics[width=3.8cm]{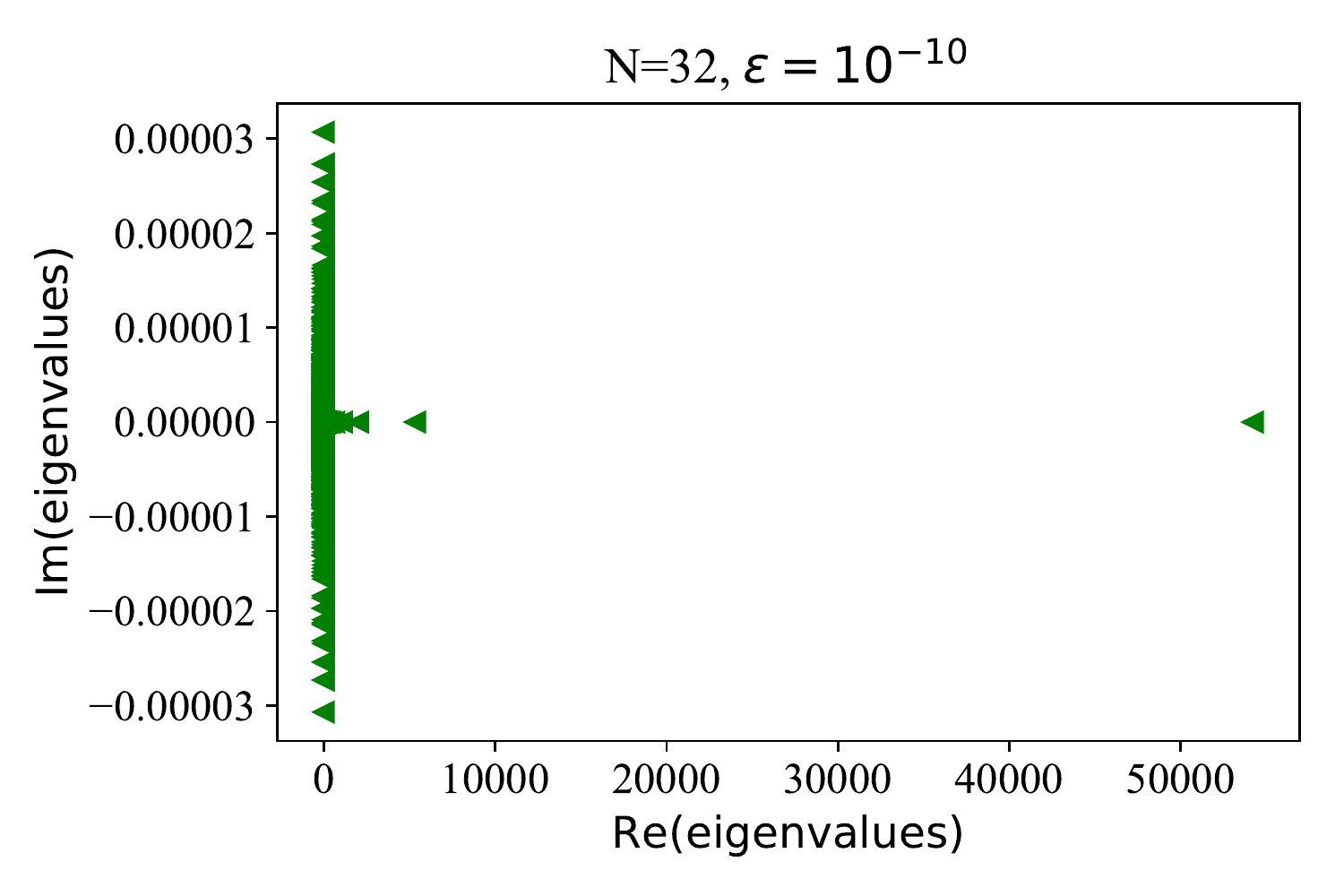} \\
\includegraphics[width=3.8cm]{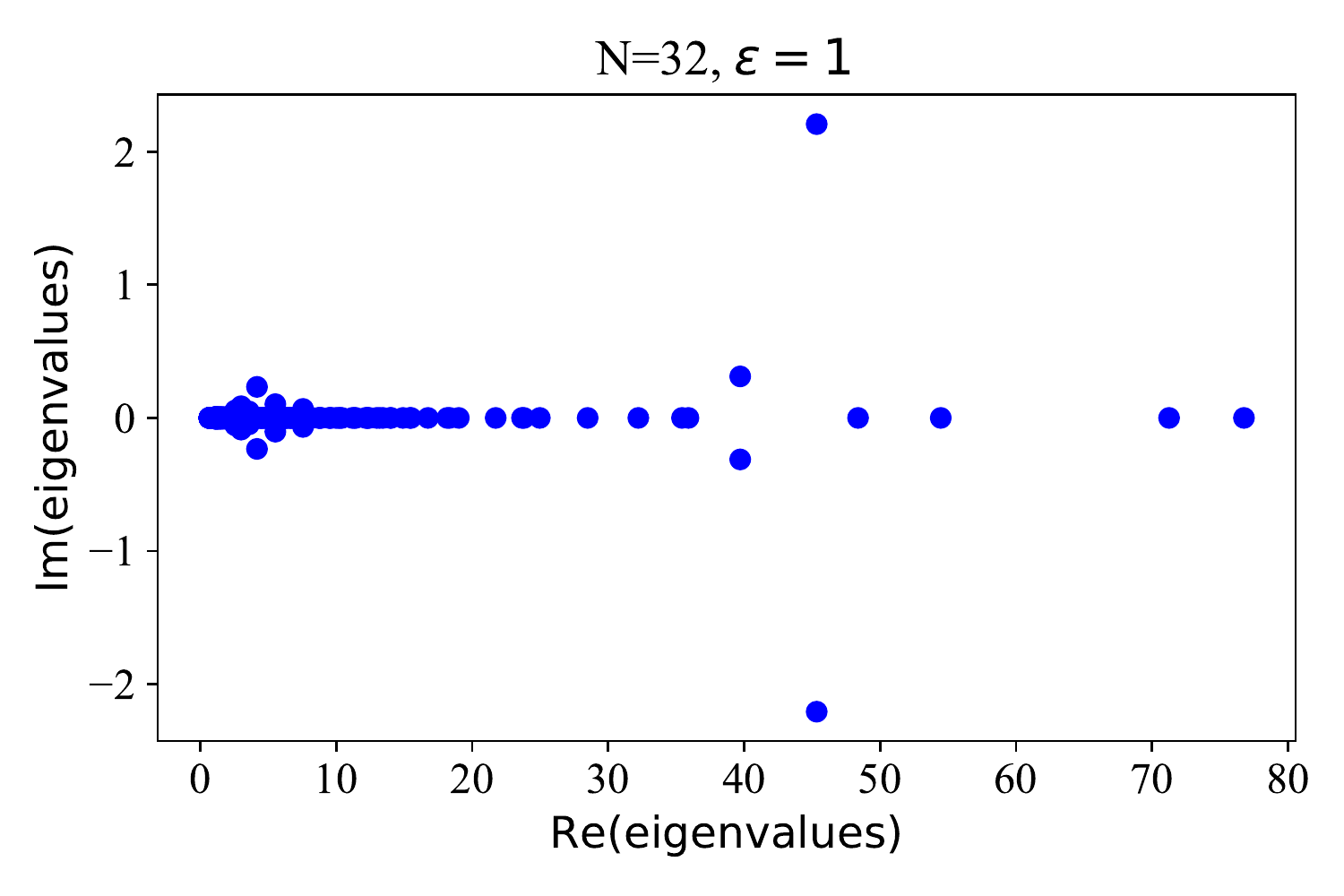} & \includegraphics[width=3.8cm]{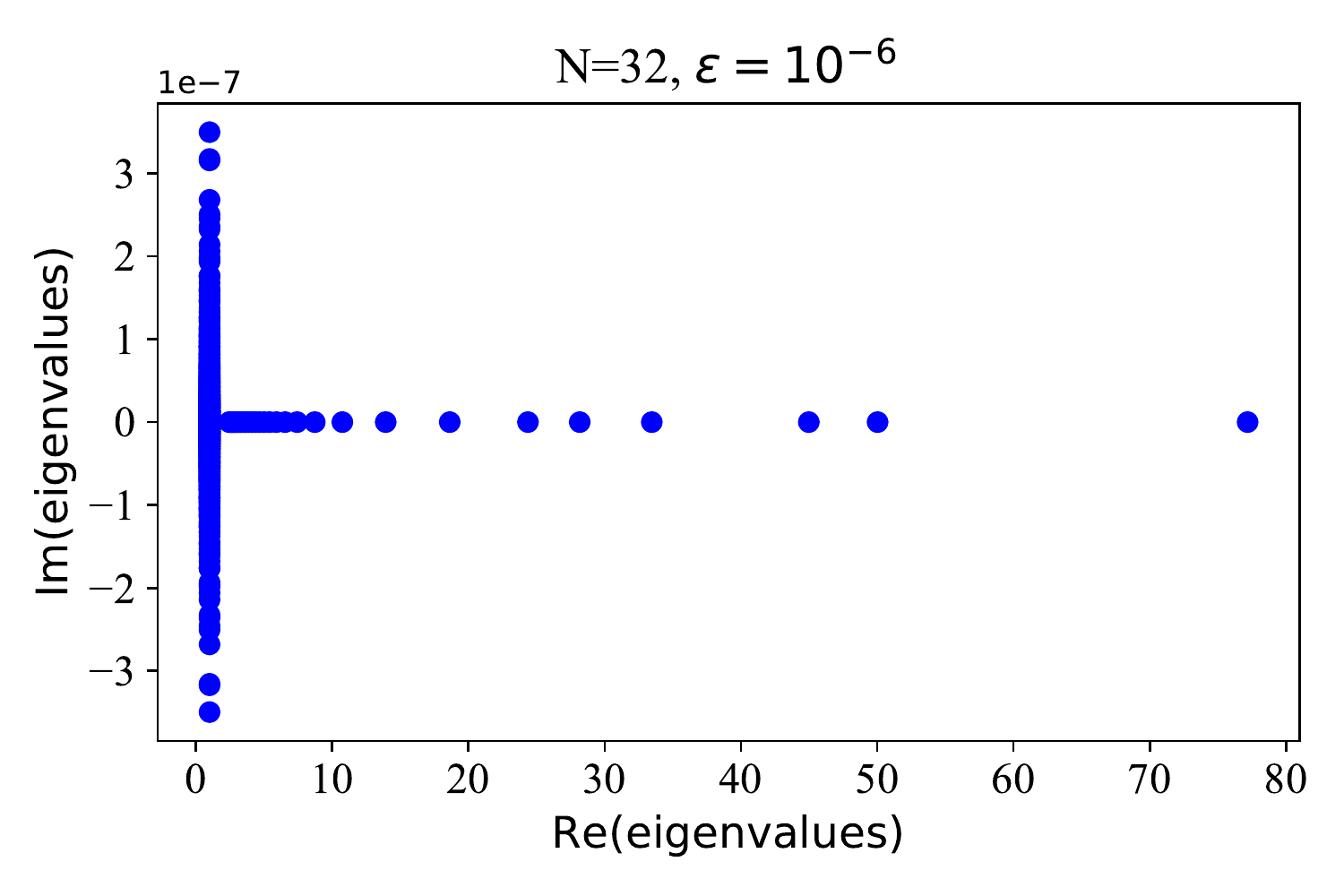} &
\includegraphics[width=3.8cm]{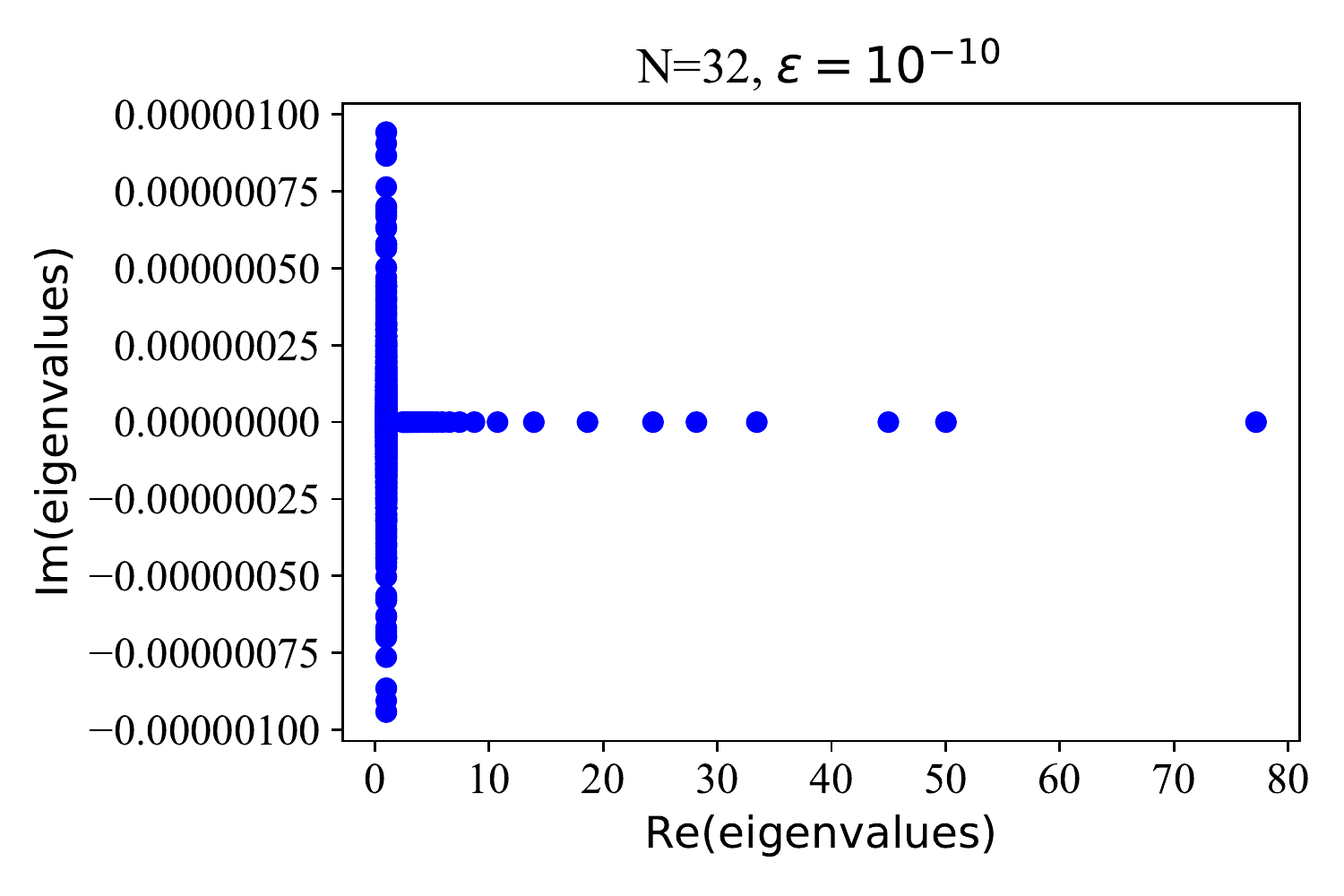} \\
\includegraphics[width=3.8cm]{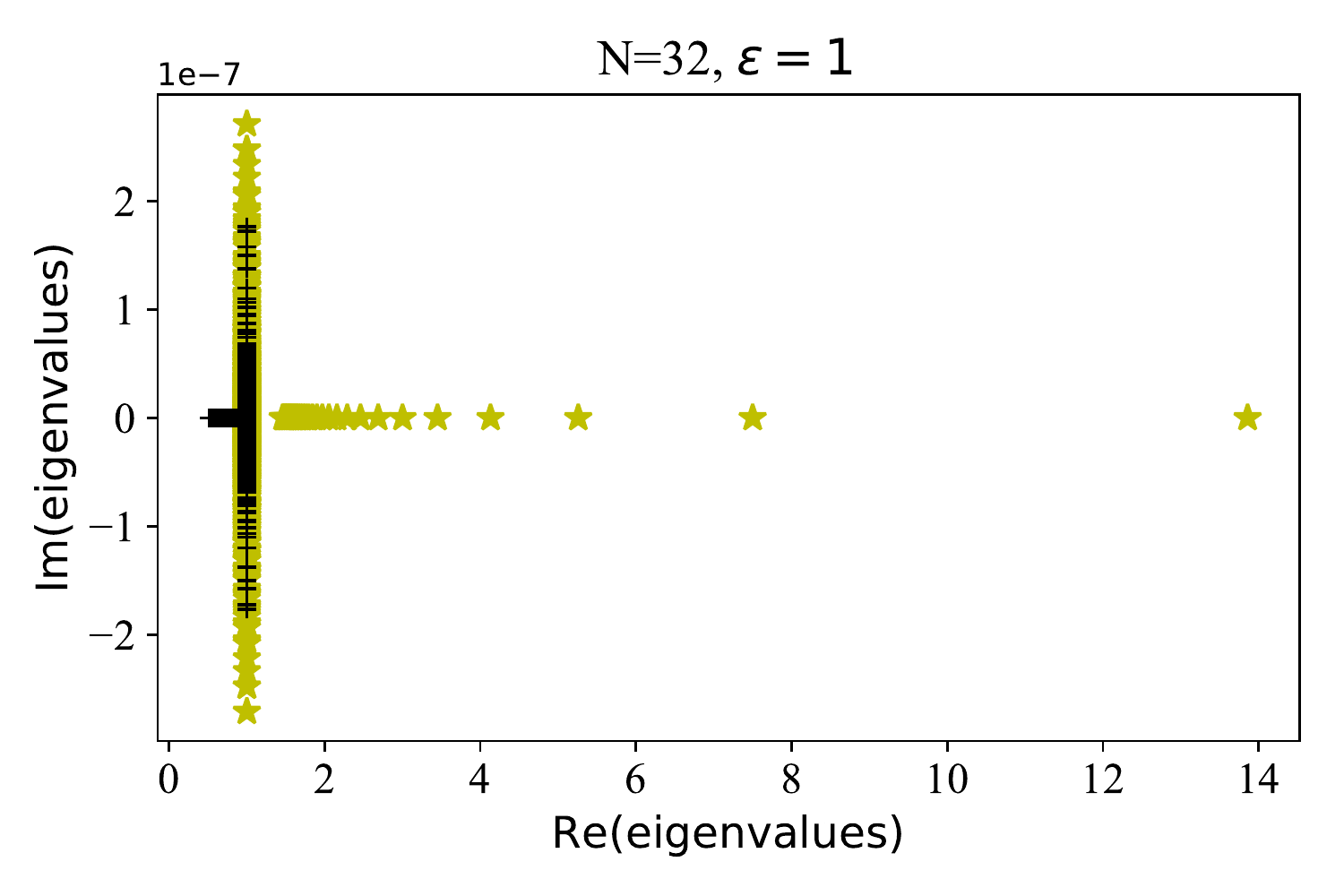} & \includegraphics[width=3.8cm]{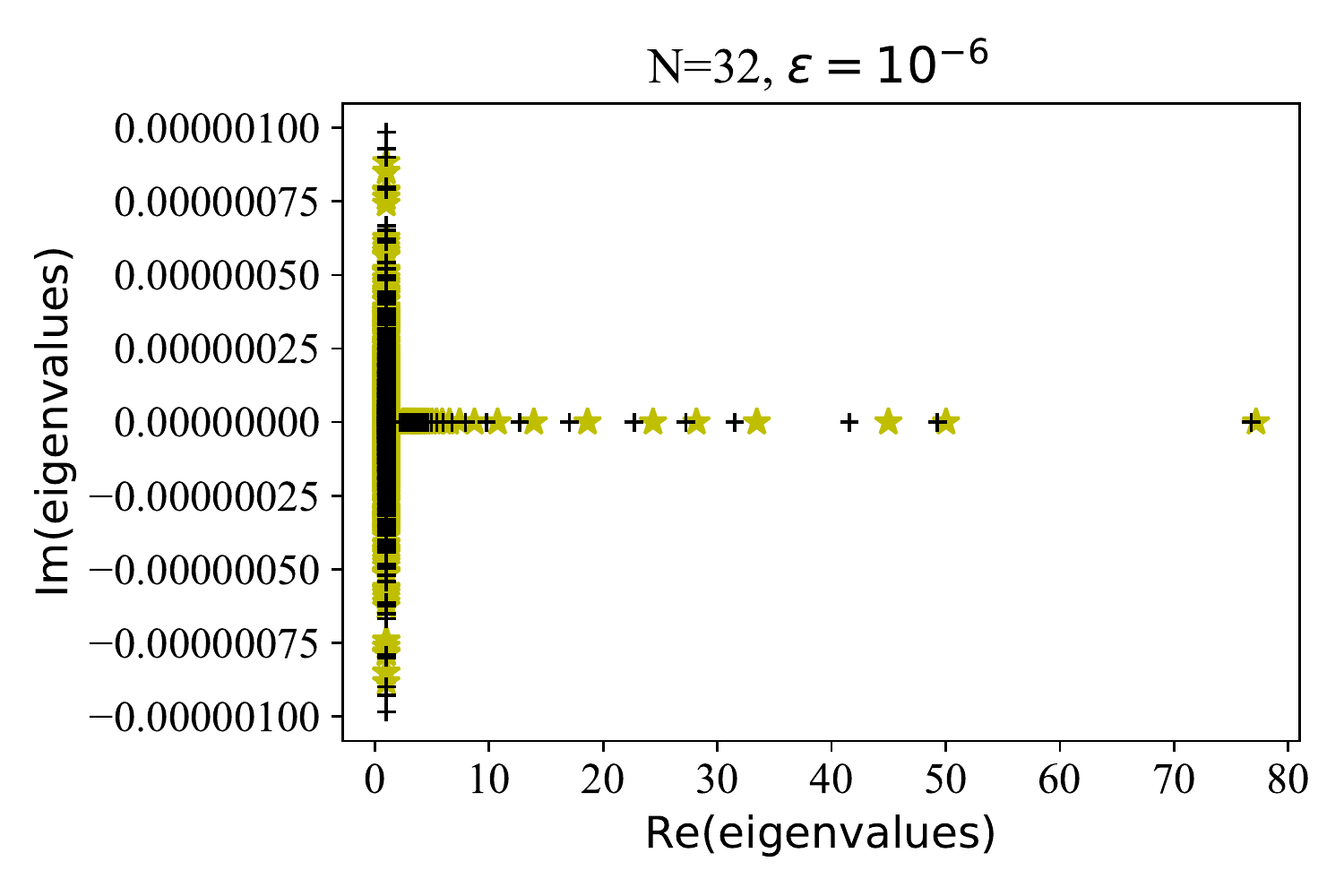} &
\includegraphics[width=3.8cm]{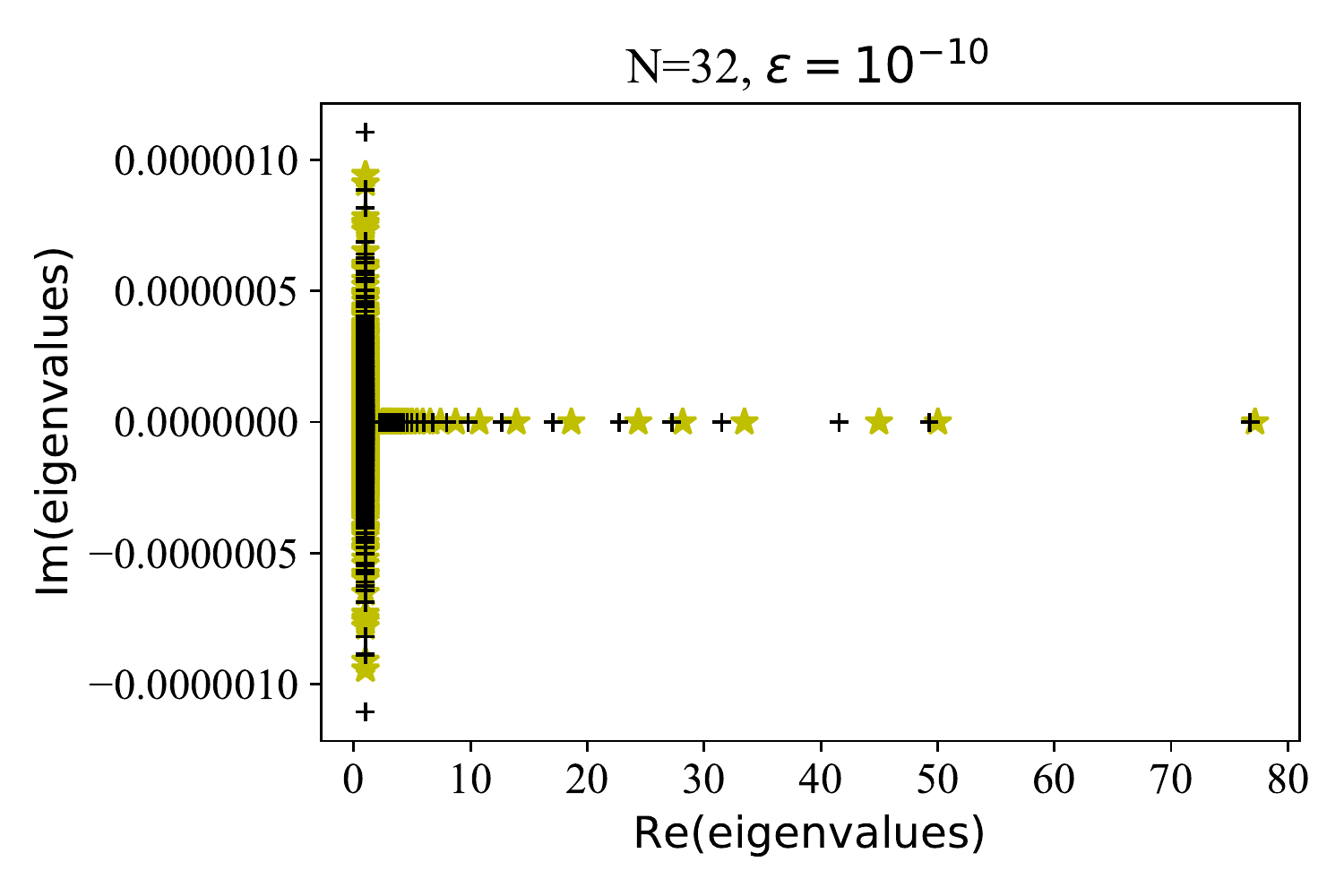}
\end{tabular}
\end{center}
 \caption{\label{fig:nonaligned_eigenvalues}
  Non-aligned case. Eigenvalues for the original system $\mathbf{A}$ and the preconditioned system $\mathbf{P}_S^{-1}\mathbf{A}$.
  Top row: unpreconditioned eigenvalues (red $\blacktriangleright$'s);
  second row: $S_1$ (green $\blacktriangleleft$'s);
  third row: $S_2$ (blue $\bullet$'s);
  bottom row: $S_3, S_5$ (yellow $\star$'s, black $+$'s respectively).
  $N$ stands for mesh size.}
\end{figure}

\begin{figure}

\begin{center}
\begin{tabular}{ccc}
\includegraphics[width=3.8cm]{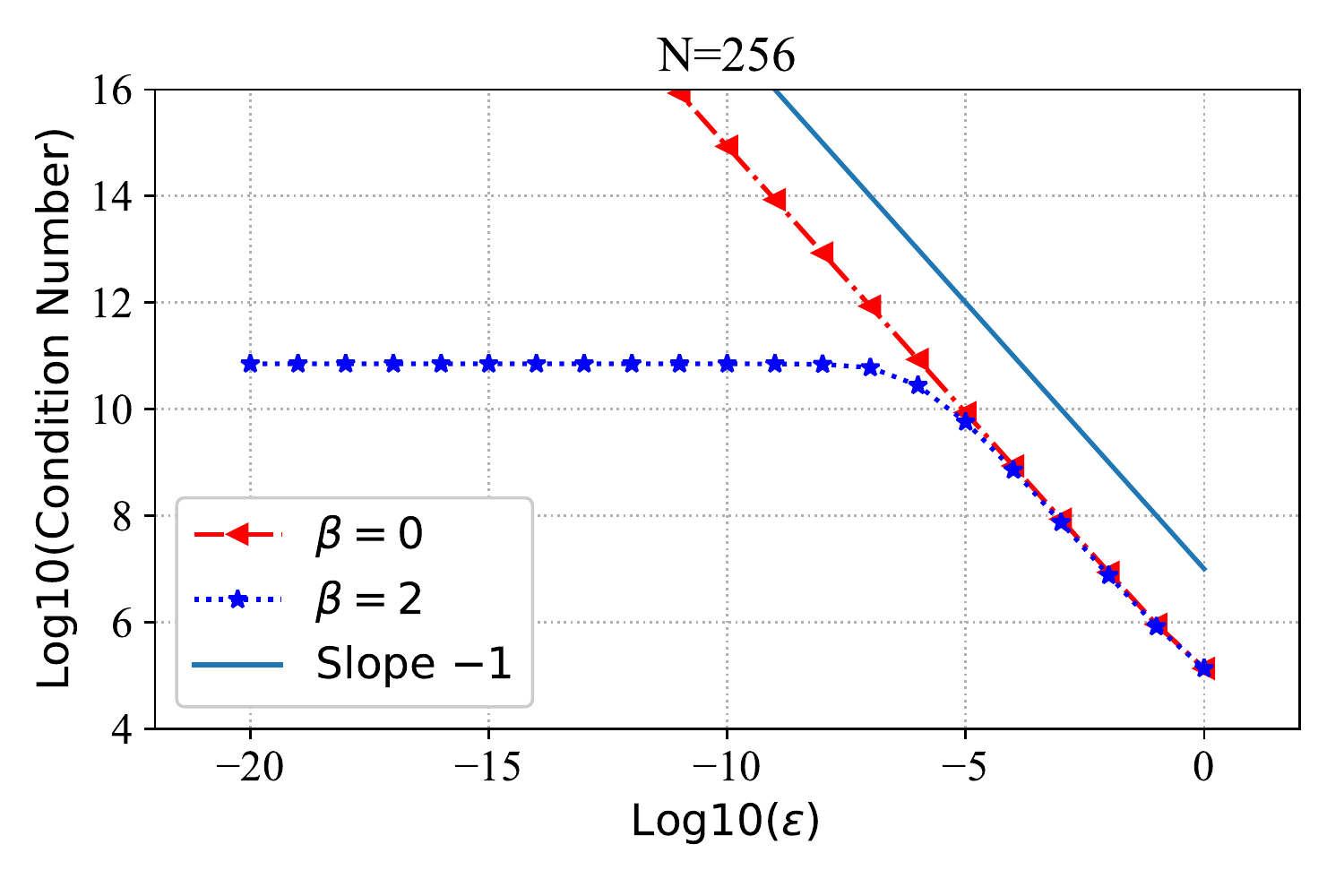} &
\includegraphics[width=3.8cm]{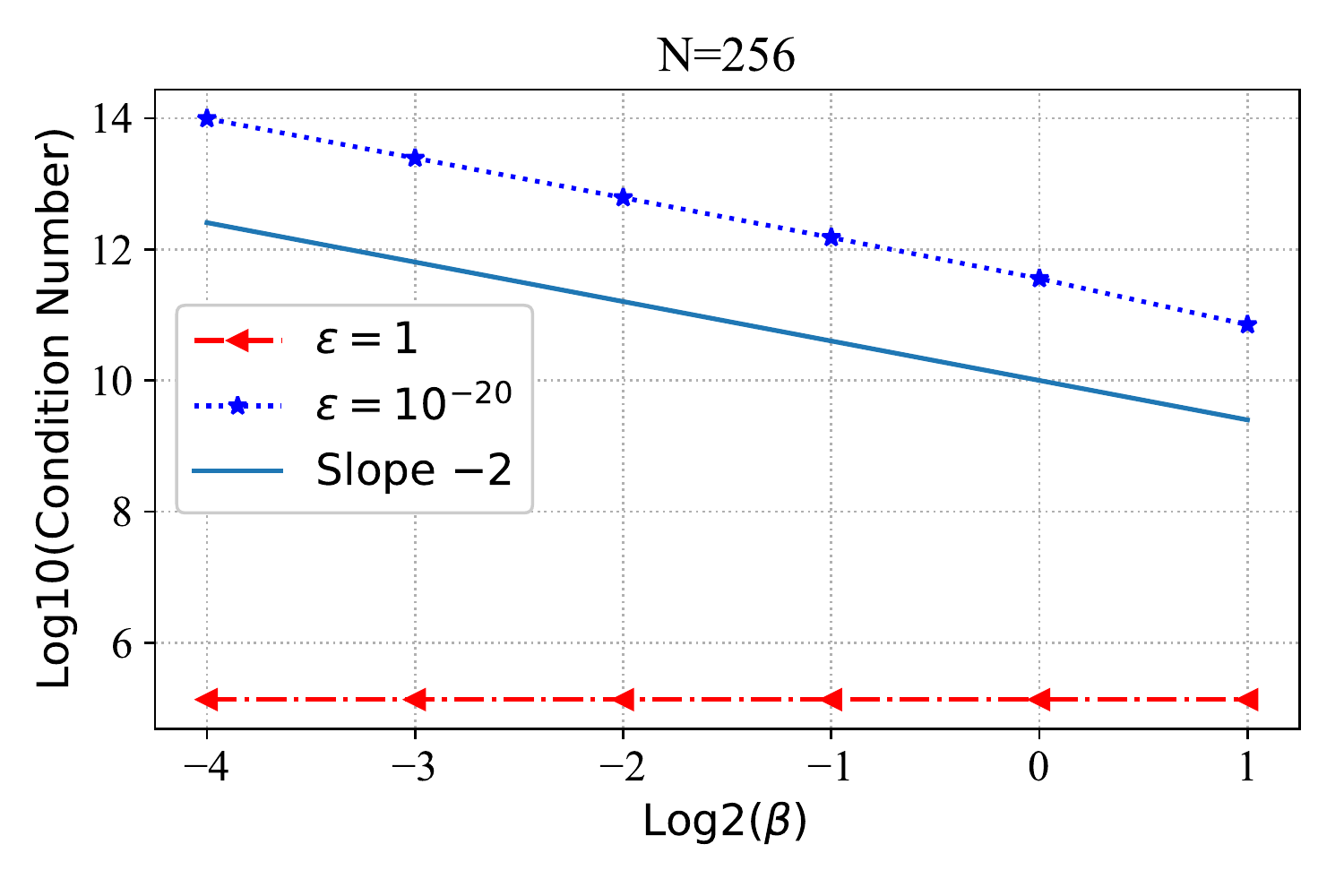} &
\includegraphics[width=3.8cm]{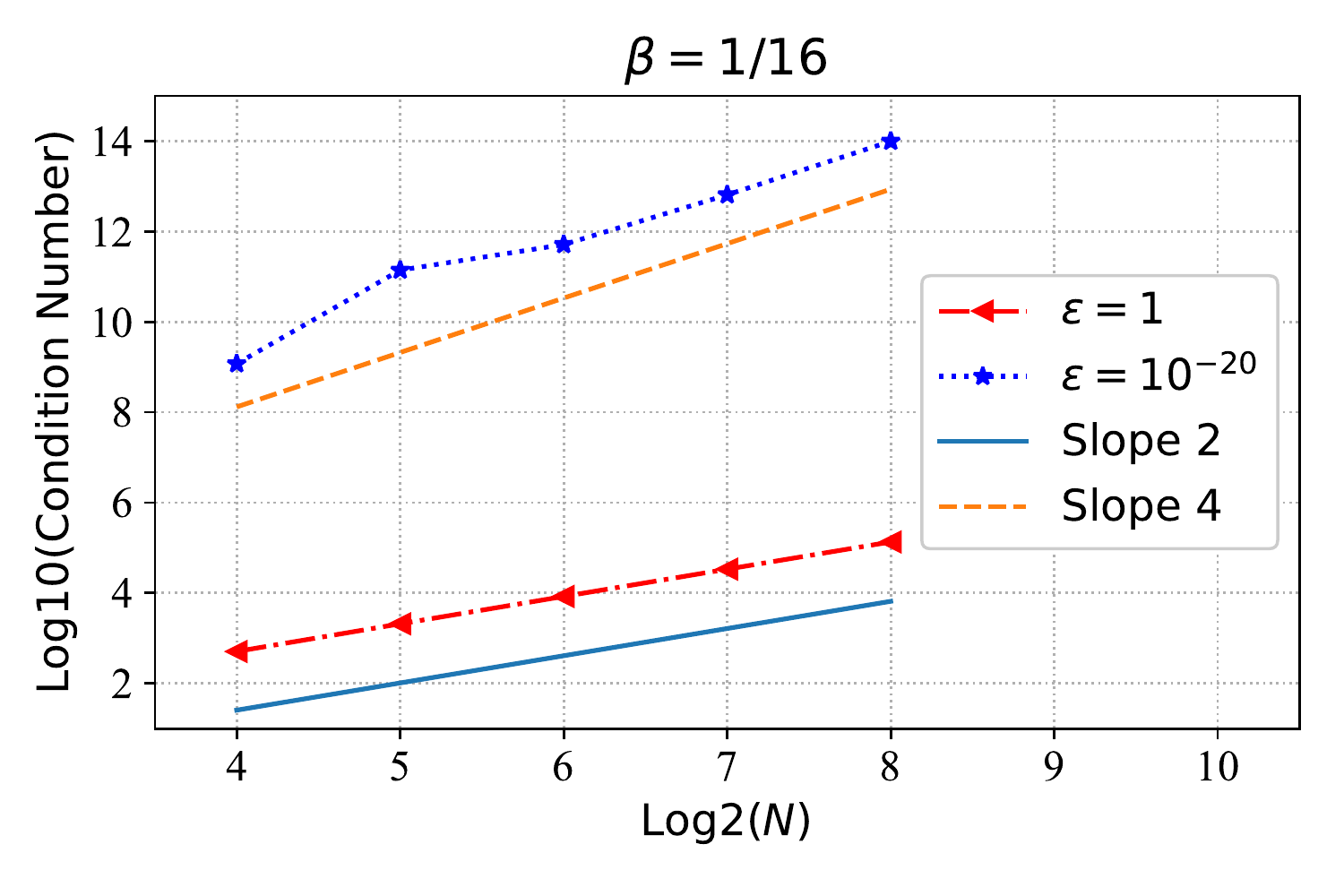} \\
(a) & (b) & (c)
\end{tabular}
\end{center}
\caption{\label{fig:nonaligned_condition}Condition number of the matrix $S_1$.
}
\end{figure}

For the non-aligned case, we use the finite difference discretization given in Appendix~\ref{sec:discrete_nonaligned}. Again, let us start with the graphics of the eigenvalues distribution of different matrices,  which is reported in Figure~\ref{fig:nonaligned_eigenvalues}. We obtain very similar graphics for the matrices $\mathbf{A}$ and   $\mathbf{P}^{-1}_{S}\mathbf{A}$ (in the cases with $S_\alpha,\alpha=2,3,5$) as we have in aligned case. One exception is for the case with $S_1$, where we find that the range of eigenvalues on the real axis is not proportional to $\frac{1}{\varepsilon}$, instead it becomes a constant  for small $\varepsilon$, {\it e.g.} $\varepsilon=10^{-10}$ or $10^{-20}$.
A very similar phenomena can also be observed from the condition number of $S_1$, which is reported in Figure~\ref{fig:nonaligned_condition}. From Figure~\ref{fig:nonaligned_condition}(a), we see that the condition number of   $S_1$ tends to be constant when $\varepsilon<10^{-6}$ for the non-aligned case, while it is proportional to $\frac{1}{\varepsilon}$ for the aligned case. Moreover, from Figure~\ref{fig:nonaligned_condition}(b), when varying $\beta$, we see that for large $\varepsilon$ the condition number is a constant  and for small $\varepsilon$ the condition number is proportional to $\frac{1}{\beta^2}$.  From Figure~\ref{fig:nonaligned_condition}(c), when varying mesh number $N$, we see that  for large $\varepsilon$ the condition number is  proportional to $N^2$  and for small $\varepsilon$ the condition number is proportional to $N^4$. Therefore, we make a conjecture of the estimation of the condition number of $S_1$ as
\begin{equation}\label{eq:condS1}
\textsf{cond}(S_1)\sim \frac{1}{\varepsilon + \beta^2 h^2}\frac{1}{h^2}.
\end{equation}
As a result, for small $\varepsilon$, the  condition number of $S_1$  is proportional to $\frac{1}{h^4}$ for fixed $\beta>0$,
thus the condition number is better than that of the aligned case. This explains that the range of real part of eigenvalues of $\mathbf{P}_S^{-1}\mathbf{A}$ does not explode  for the non-aligned case when vanishing $\varepsilon$ and
the good performance even with small $\varepsilon$. Finally, from numerical computation, we find that $S_4$ is almost equal to the exact Schur complement $E$, thus the corresponding eigenvalues are very tightly clustered to 1. This fact indicates that the difference between the exact Schur complement $E$ and $A_1+\varepsilon B$ is just the added rows corresponding to the boundary $[0,1]\times\{0\}$, therefore the more
accurate we approximate well these boundary rows, the better performance of the block preconditioner $\mathbf{P}_S$ we can expect.

The numbers of the GMRES iteration are summarized in Table~\ref{tab:gmresiternum_nonaligned}. Remarkably for the case with $S_1$, we see 
the GMRES iterations almost do not increase when refining meshes and they increase slightly when $\varepsilon$ vanishing. It seems that the case with  $S_1$ preforms the best comparing to the cases with $S_\alpha,\alpha=2,3,5$.
This phenomena is very different from the aligned case, especially when $\varepsilon=10^{-20}$, where the GMRES method does not converge.
This "robustness" seems to verify our conjecture \eqref{eq:condS1} from a different point of view.
For the case with $S_2$, we see again the GMRES method does not work well for large $\varepsilon$.
For $S_3$, the GMRES iteration are almost independent of the choice of $\varepsilon$.
More importantly, it does not deteriorate again for vanishing $\beta$.
For the case with $S_4$, the GMRES method converges in two iterations which indicates that $S_4$ is almost equal to the exact Schur complement $E$.
Finally, the case of $S_5$ combines the advantages of the cases with $S_1$ and $S_3$.
Thus we conclude that approximate Schur complement $S_1$ and $S_5$ are the preferred ones for non-aligned case.

\begin{table}[!h]
\caption{\label{tab:gmresiternum_nonaligned}
  Non-aligned case. Iteration number of GMRES method with approximate Schur complements $S_\alpha$, $\alpha=1,\dots,5$.}
\begin{center}

(a)  Case with $S_1$

\begin{tabular}{ccccccc}
\hline
$\varepsilon$ & 1 & $10^{-1}$ & $10^{-2}$ & $10^{-6}$ & $10^{-10}$ & $10^{-20}$\\
\hline
$N = 32$ & 5 & 8 & 14 & 13 & 13 & 13 \\
$N = 64$ & 5 & 7 & 11 & 15 & 15 & 15 \\
$N = 128$ & 4 & 6 & 9 & 17 & 17 & 17 \\
$N = 256$ & 4 & 5 & 7 & 19 & 18 & 19 \\
\hline
\end{tabular}

\vspace{3mm}

(b)  Case with $S_2$

\begin{tabular}{ccccccc}
\hline
$\varepsilon$ & 1 & $10^{-1}$ & $10^{-2}$ & $10^{-6}$ & $10^{-10}$ & $10^{-20}$\\
\hline
$N = 32$ & 33 & 22 & 18 & 12 & 12 & 12 \\
$N = 64$ & 50 & 33 & 25 & 15 & 15 & 15 \\
$N = 128$ & 73 & 48 & 34 & 17 & 17 & 17 \\
$N = 256$ & $\times$ & 66 & 46 & 19 & 20 & 20 \\
\hline
\end{tabular}

\vspace{3mm}
(c)  Case with $S_3$

\begin{tabular}{ccccccc}
\hline
$\varepsilon$ & 1 & $10^{-1}$ & $10^{-2}$ & $10^{-6}$ & $10^{-10}$ & $10^{-20}$\\
\hline
$N = 32$ & 11 & 12 & 14 & 13 & 12 & 12 \\
$N = 64$ & 13 & 15 & 17 & 15 & 15 & 15 \\
$N = 128$ & 15 & 18 & 21 & 20 & 17 & 17 \\
$N = 256$ & 17 & 21 & 25 & 25 & 20 & 20 \\
\hline
\end{tabular}

\vspace{3mm}

(d)  Case with $S_4$

\begin{tabular}{ccccccc}
\hline
$\varepsilon$ & 1 & $10^{-1}$ & $10^{-2}$ & $10^{-6}$ & $10^{-10}$ & $10^{-20}$\\
\hline
$N = 32$ & 2 & 2 & 2 & 2 & 2 & 2 \\
$N = 64$ & 2 & 2 & 2 & 2 & 2 & 2 \\
$N = 128$ & 2 & 2 & 2 & 2 & 2 & 2 \\
$N = 256$ & 2 & 2 & 2 & 2 & 2 & 2 \\
\hline
\end{tabular}

\vspace{3mm}
(e) Case with $S_5$

\begin{tabular}{ccccccc}
\hline
$\varepsilon$ & 1 & $10^{-1}$ & $10^{-2}$ & $10^{-6}$ & $10^{-10}$ & $10^{-20}$\\
\hline
$N = 32$ & 7 & 11 & 13 & 12 & 12 & 12 \\
$N = 64$ & 6 & 12 & 16 & 15 & 14 & 14 \\
$N = 128$ & 6 & 12 & 18 & 19 & 17 & 17 \\
$N = 256$ & 6 & 12 & 21 & 25 & 19 & 19 \\
\hline
\end{tabular}

\end{center}
\end{table}

\begin{figure}

\begin{center}
\begin{tabular}{cc}
\includegraphics[width=5cm]{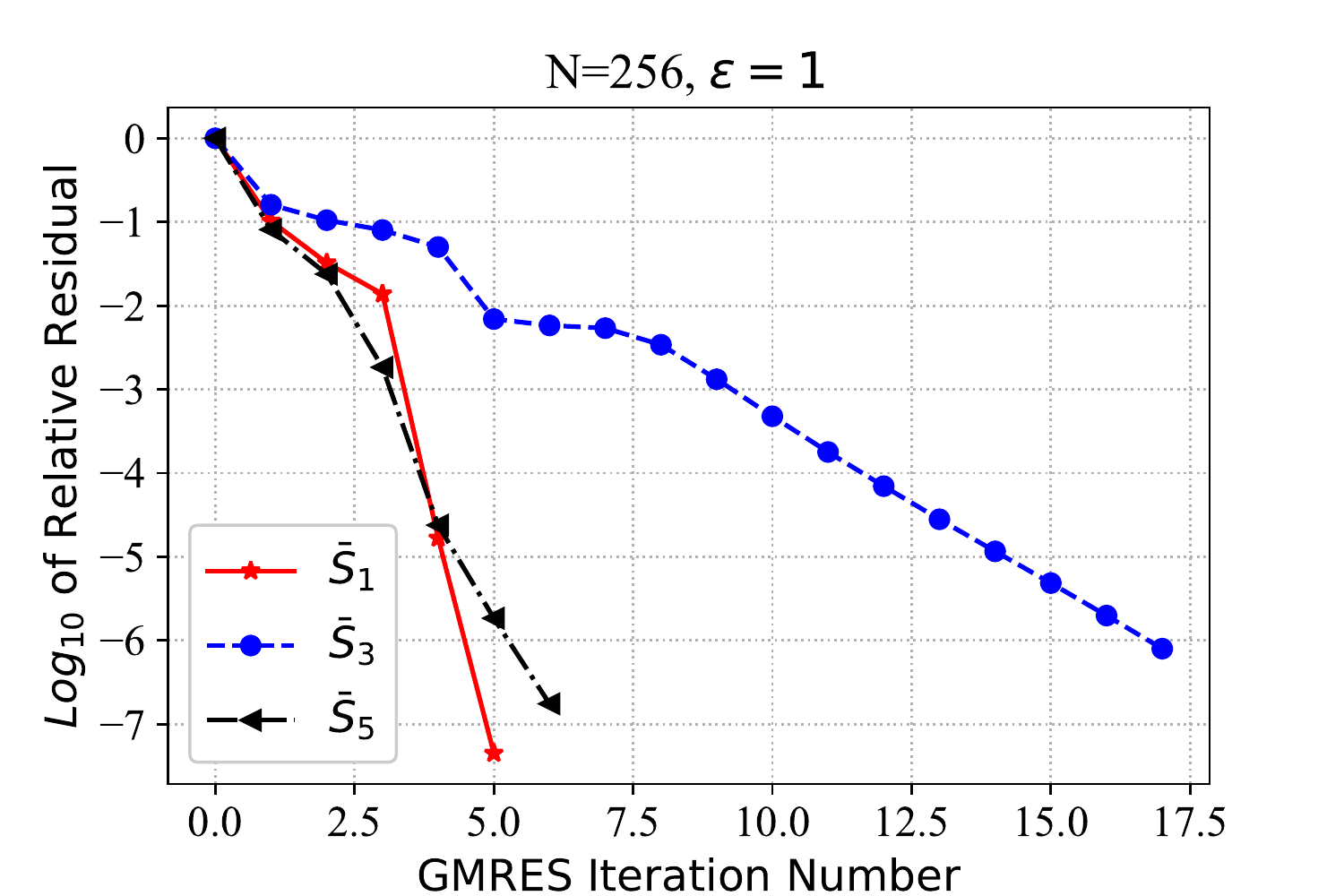} & \includegraphics[width=5cm]{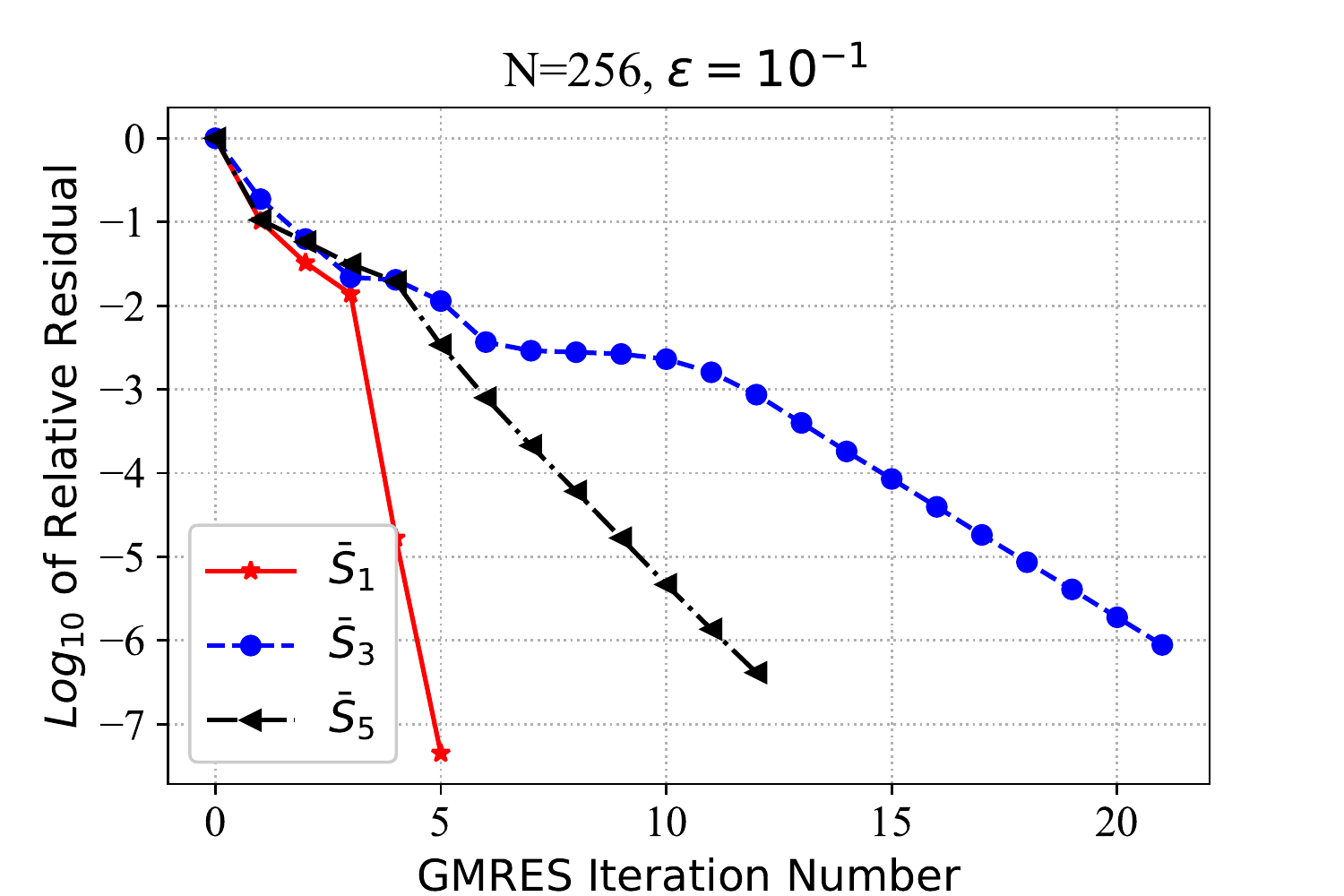} \\
\includegraphics[width=5cm]{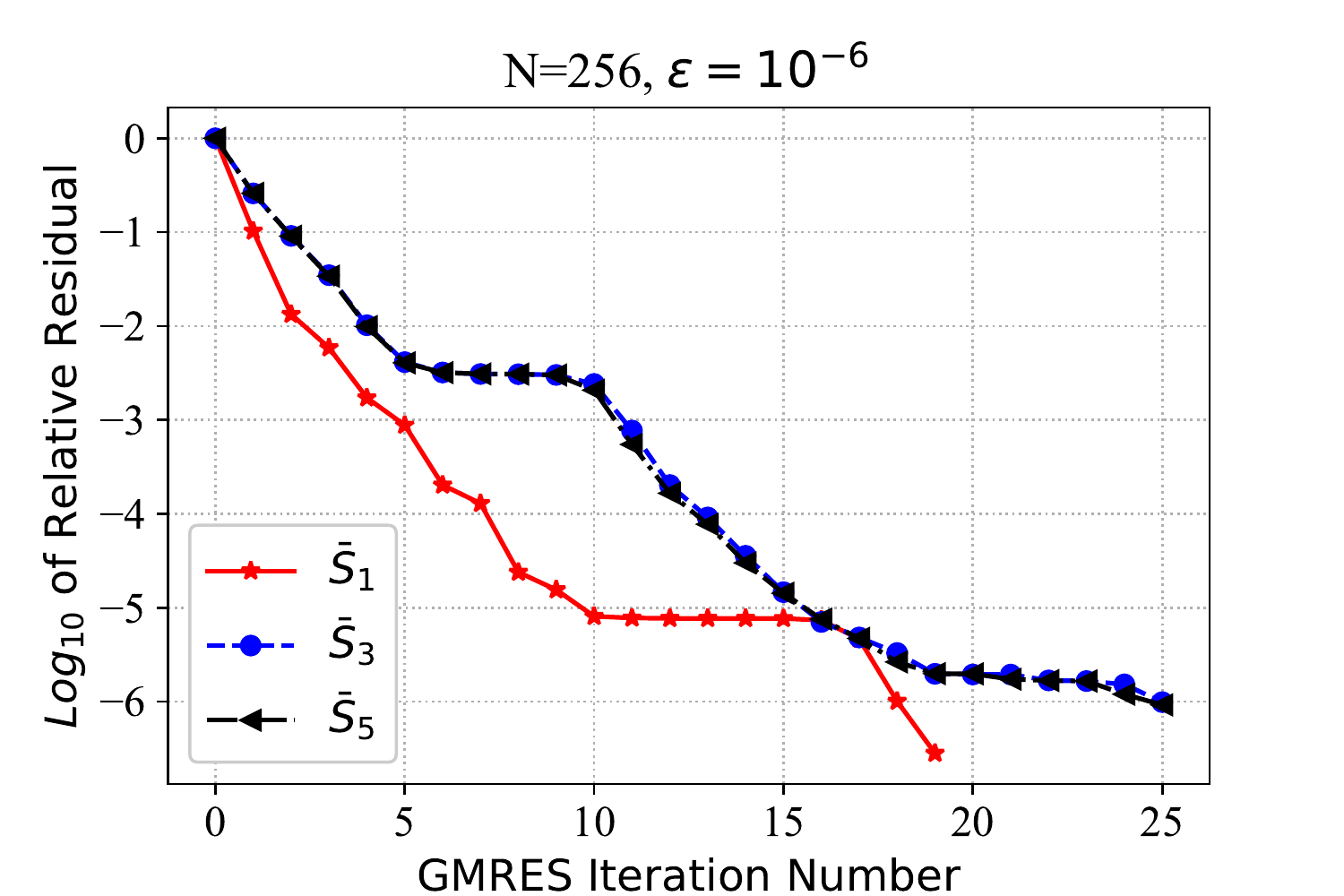} & \includegraphics[width=5cm]{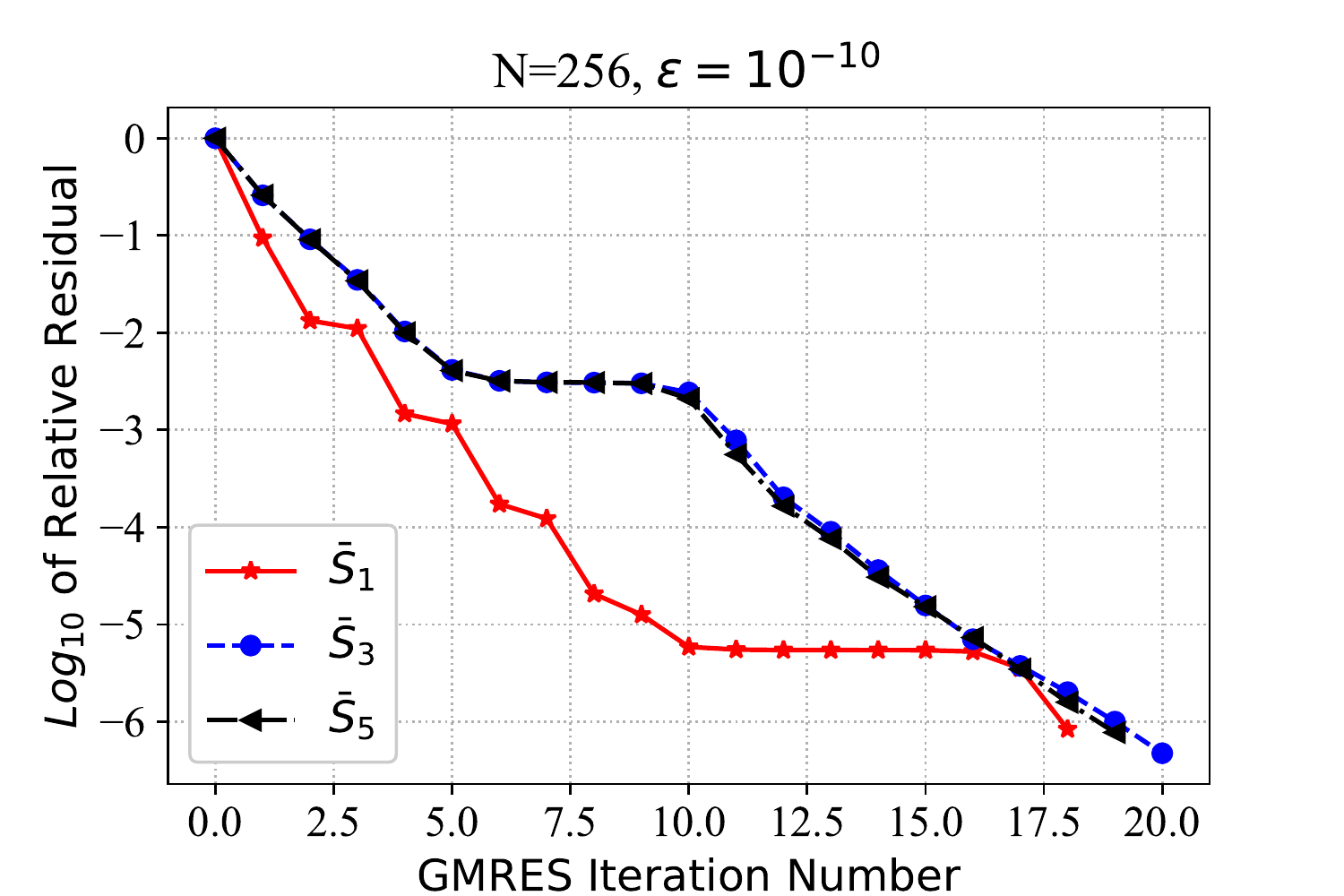}
\end{tabular}
\end{center}
\caption{\label{fig:nonaligned_res}Convergence history of GMRES iteration for the non-aligned case.}
\end{figure}

Figure~\ref{fig:nonaligned_res} presents relative residual history of GMRES iteration. We see a clear tendency of convergence of GMRES method with different approximate Schur complements. For the current case with $\beta=2$, we find $S_1$ performs the best. However, from the conjecture~\eqref{eq:condS1}, the condition number of $S_1$ deteriorates when vanishing $\beta$, thus the convergence of the GMRES method may also be perturbed. For instance, by setting $\beta=0.1$, the GMRES method does not converge when $N=256$.
The case with $S_5$ performs less well than the one with $S_1$, but it is better than the case with $S_3$, especially for large $\varepsilon$.

\begin{figure}

\begin{center}
\begin{tabular}{cc}
\includegraphics[width=5cm]{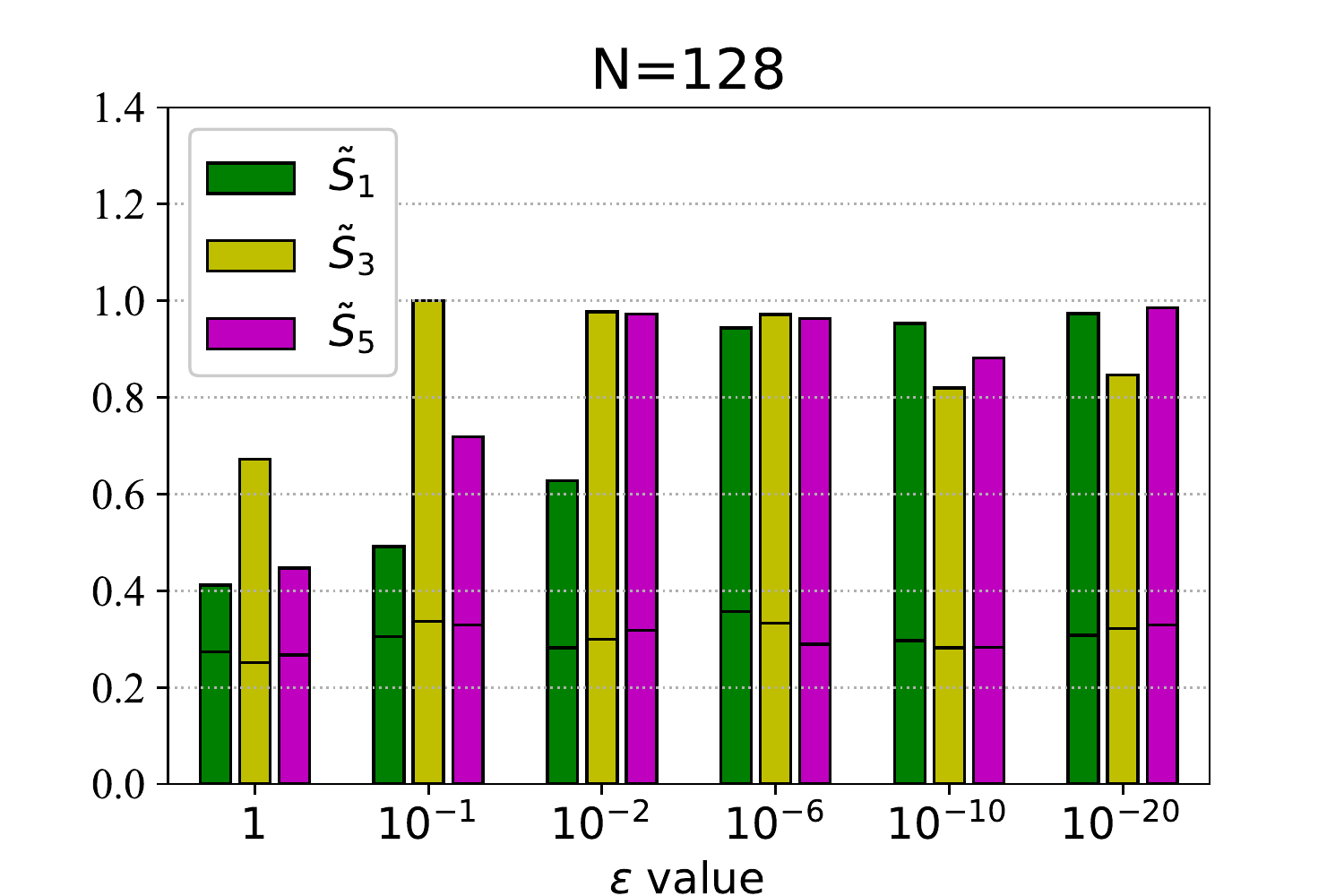} & \includegraphics[width=5cm]{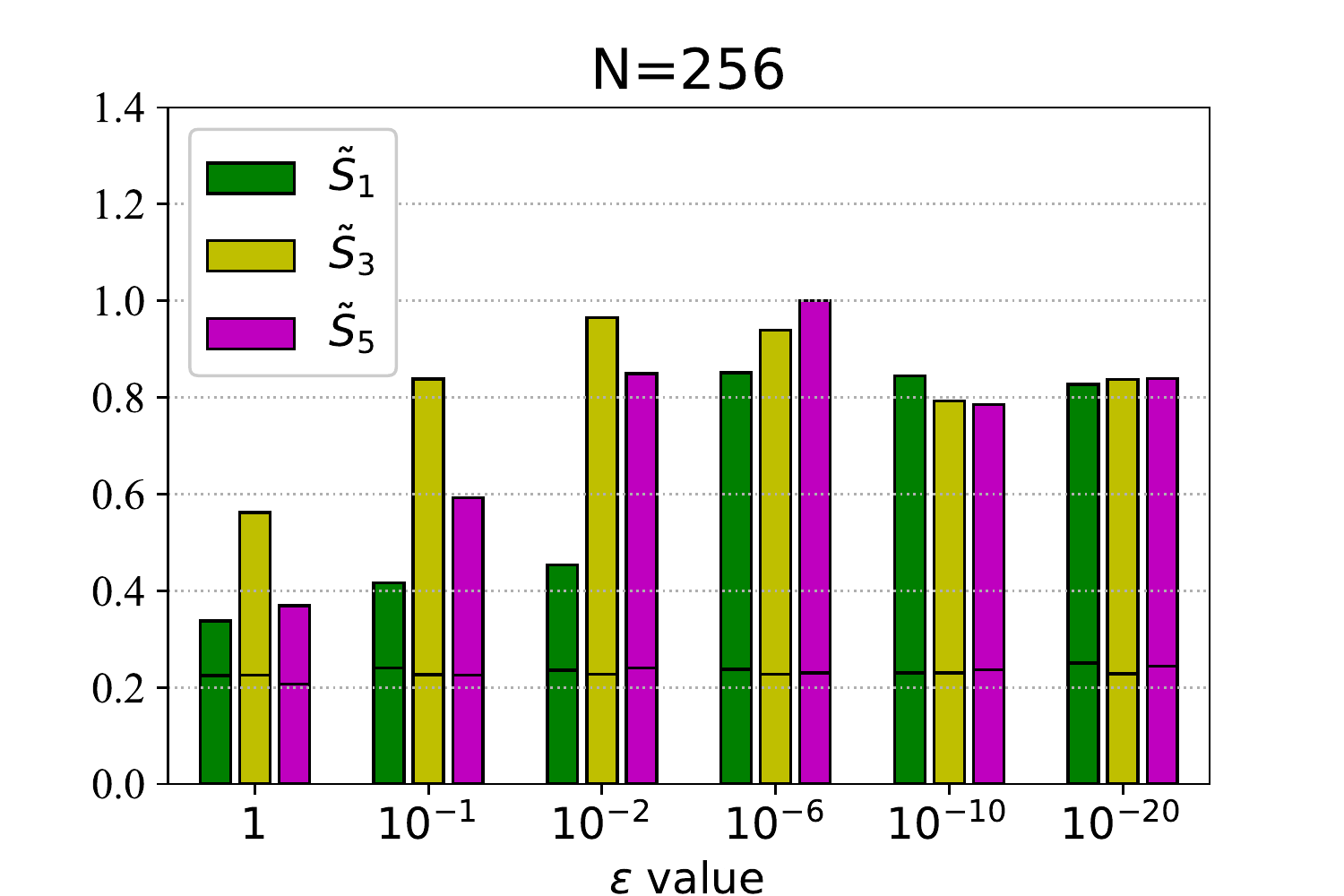}
\end{tabular}
\end{center}
\caption{\label{fig:nonaligned_directtime}
  Non-aligned case. The computational time for solving preconditioned linear system. The presented time is rescaled by the maximum value from each group of comparison. $\tilde{S}_\alpha,\alpha=1,3,5$ stands for different approximate Schur complement. Each bar is decomposed into two portions: bottom one is time of matrix factorizations and top one is for resolution time. $N$ stands for mesh size.}
\end{figure}

The computational time for solving preconditioned linear system is reported in Figure~\ref{fig:nonaligned_directtime}. Again, we factorize first the matrices $A_3$ and $S_\alpha$, then we apply the GMRES method. The first observation is that, very different from the aligned case, for the non-aligned case most computational time comes from the GMRES solver. For large $\varepsilon$, the cases with $S_1$ and $S_5$ are more efficient than the one with $S_3$. But for small $\varepsilon$, the efficiency of all choices is almost the same.
Therefore we recommend to use the preconditioner with the approximate Schur complement $S_5$.

\subsection{Without inflow boundary condition}
In Section~\ref{sec:MMAPscheme}, the introduced MMAP scheme provides a unique couple of $(\phi,q)$, where the key point is to impose the inflow condition for $q$. However, we notice that the auxiliary  $q$ does not involve directly in computation of $\phi$, instead  the parallel derivative $\nabla_\| q$ is used in the MMAP scheme.
This similar situation also appears in the solution of saddle point problem from incompressible
Navier-Stokes equations (see~\cite{nsref1} and references therein), where the matrix of NS equations is singular because pressure variable
is not unique even though its gradient is unique.

In preconditioned GMRES iteration, one mainly needs to apply matrix vector multiplication
and a inverse of the preconditioner.
Thus in principle GMRES method can still work even if the matrix is singular provided
the preconditioner is applicable. 
In the solving of NS equations~\cite{nsref1}, GMRES can obtain right solutions for the velocity.
In this paper $q$ is just a auxiliary variable, thus the inflow condition can be dropped if GMRES iteration can generate
the correct solutions for $\phi$. 
In another word, the inflow condition is not necessary to be introduced for determining a unique $\phi$ when using GMRES. But when using direct solvers we can not drop the inflow condition.

In this subsection, we substitute the inflow condition in~\eqref{eq:MM1} by $\mathbf{n}\cdot \nabla_\| q + \mathbf{n}\cdot \nabla_\bot \phi =0,(x,z)\in [0,1]\times\{0\}$ and solve
\begin{equation}\label{eq:MM3}
\left\{
\begin{array}{ll}
 - \Delta_\| q-\Delta_\bot\phi = f,& (x,z)\in(0,1)\times(0,1),\\[3mm]
q=0,&(x,z)\in\{0,1\}\times[0,1],\\[3mm]
\mathbf{n}\cdot \nabla_\| q + \mathbf{n}\cdot \nabla_\bot \phi =0,&(x,z)\in [0,1]\times\{0,1\}.
\end{array}
\right.
\end{equation}
Now the MMAP scheme becomes~\eqref{eq:MM3} and \eqref{eq:MM2}, and its corresponding discretization matrix (denoted again by $\mathbf{A}$) is no longer invertible, because the discretization matrix of the operator $-\Delta_\|q$ of~\eqref{eq:MM3} (denoted by  $\hat{A}_3$) is now not invertible. 
We wirte the new linear system by
\begin{equation}\label{eq:MMAPlinsyst2}
\begin{pmatrix}
\hat{A}_3 & B \\
A_2 & A_1
\end{pmatrix}
\begin{pmatrix}
\hat{Q} \\
\Phi
\end{pmatrix}
=
\begin{pmatrix}
F \\
0
\end{pmatrix}.
\end{equation}
Remember that we need the preconditioner to be applicable.
To make the sub-blocks of the preconditioner invertible,
we use $A_3$ to compute the approximate Schur complement.
Namely we use $A_3$ as an approximation of $\hat{A}_3$, so that the preconditioners discussed in Section~\ref{sec:blockprecon} remain the same.
In short we use inflow conditions in the preconditioning but not in the computation of coefficient matrix $\mathbf{A}$.

Here, we only consider the approximate Schur complement $S_5$, which is shown to be the most robust choice in previous subsections. 
The iteration number of the GMRES method is reported in Table~\ref{tab:gmresiternum_noinflow}. For aligned case, Table~\ref{tab:gmresiternum_noinflow}(a) is the same as in Table~\ref{tab:gmresiternum_aligned}(e). While for non-aligned case, Table~\ref{tab:gmresiternum_noinflow}(b) is also similar as in Table~\ref{tab:gmresiternum_nonaligned}(e).
This show that even though the matrix of linear system~\eqref{eq:MMAPlinsyst2} is not invertible,
the proposed preconditioned GMRES algorithm still works well.

\begin{table}[!h]
\caption{\label{tab:gmresiternum_noinflow}
Iteration number of the GMRES method for the scheme without inflow boundary condition. The approximate Schur complement $S_5$ is used.}
\begin{center}

(a) Aligned case

\begin{tabular}{ccccccc}
\hline
$\varepsilon$ & 1 & $10^{-1}$ & $10^{-2}$ & $10^{-6}$ & $10^{-10}$ & $10^{-20}$\\
\hline
$N = 32$ & 4 & 4 & 4 & 3 & 2 & 2 \\
$N = 64$ & 4 & 4 & 4 & 3 & 2 & 2 \\
$N = 128$ & 4 & 4 & 4 & 3 & 2 & 2 \\
$N = 256$ & 4 & 4 & 4 & 4 & 2 & 2 \\
\hline
\end{tabular}

\vspace{3mm}
(b) Non-aligned case

\begin{tabular}{ccccccc}
\hline
$\varepsilon$ & 1 & $10^{-1}$ & $10^{-2}$ & $10^{-6}$ & $10^{-10}$ & $10^{-20}$\\
\hline
$N = 32$ & 6 & 10 & 13 & 12 & 12 & 12 \\
$N = 64$ & 6 & 12 & 16 & 15 & 14 & 14 \\
$N = 128$ & 6 & 12 & 18 & 19 & 17 & 17 \\
$N = 256$ & 6 & 12 & 20 & 25 & 19 & 19 \\
\hline
\end{tabular}

\end{center}
\end{table}

\section{Conclusion and perspectives}
In this paper, we introduce  block preconditioning methods for the MMAP method.
We propose a series of block preconditioners based on the approximate Schur complement techniques.
We perform a number of numerical tests for both the aligned case and non-aligned case.
Numerical results show that on the one hand, the proposed block preconditioners are robust with grid
refinement; on the other hand, they can work well for all anisotropic strengths.

The MMAP method works for the cases where the border of magnetic field lines cross the boundaries of computational domain.
For the case with closed magnetic field lines, we still need to explore other asymptotic-preserving methods.
For instance, in~\cite{2021Li}, an asymptotic-preserving scheme based on first order system (APFOS) was introduced.
Therefore, to develop a suitable block preconditioning method for APFOS scheme can be an extension of present work.
Moreover the implementation of finite element method and better approximation of the boundary rows of
exact Schur complement is under consideration.
\section*{Acknowledgements}
This work has been supported by Heilongjiang Provincial Natural Science Foundation of China (LH2019A013) and the Fundamental Research Funds for the Central Universities.

Chang Yang acknowledges support by a public grant from the "Laboratoire d'Excellence Centre International de Math\'ematiques et d'Informatique" (LabEx CIMI) overseen by the French National Research Agency (ANR) as part of the "Investissements d'Avenir" program (reference ANR-11-LABX-0040) in the frame of the PROMETEUS project (PRospect of nOvel nuMerical modEls for elecTric propulsion and low tEmperatUre plaSmas) and from the FrFCM (F\'ed\'eration de recherche pour la Fusion par Confinement Magn\'etique) in the frame of the NEMESIA project (Numerical mEthods for Macroscopic models of magnEtized plaSmas and related anIsotropic equAtions).

Lingxiao Li has also been supported by National Natural Science Foundation of China under Grant 11901042.

%
  \appendix
\section{Finite difference discretizations}\label{sec:FDD}
In this part, let us give the finite difference discretizations for the MMAP scheme introduced in Section~\ref{sec:MMAPscheme}.

\subsection{Aligned case}\label{sec:discrete_aligned}
In the aligned case, we consider the magnetic field aligned with $z$-axis, {\it i.e.} $\mathbf{b}=(0,1)^T$. For the convenience of construction of matrices, we write the MMAP scheme of Section~\ref{sec:MMAPscheme} as follows
\begin{equation}\label{eq:MM1aligned}
\left\{
\begin{array}{ll}
 - \partial_z^2 q - \partial_x^2\phi = f,& (x,z)\in(0,1)\times(0,1),\\[3mm]
q=0,&(x,z)\in\{0,1\}\times[0,1],\\[3mm]
q=0,&(x,z)\in [0,1]\times\{0\}, \\[3mm]
\partial_z q=0,&(x,z)\in [0,1]\times\{1\},
\end{array}
\right.
\end{equation}
and
\begin{equation}\label{eq:MM2aligned}
\left\{
\begin{array}{ll}
\varepsilon\partial_z^2 q - \partial_z^2\phi   = 0,& (x,z)\in(0,1)\times(0,1),\\[3mm]
\phi=0,&(x,z)\in\{0,1\}\times[0,1],\\[3mm]
\partial_z q=0,&(x,z)\in [0,1]\times\{0\},\\[3mm]
\partial_z\phi = 0,& (x,z)\in[0,1]\times\{1\}.
\end{array}
\right.
\end{equation}

The mesh is defined by
\begin{equation*}
\left\{
\begin{array}{ll}
x_i = (i-1)\Delta x,\ i\in\{1,\dots,n_x\},\\[3mm]
z_j = (j-3/2)\Delta z,\ j\in\{1,\dots,n_z\},
\end{array}
\right.
\end{equation*}
where
$\Delta x = \frac{1}{n_x-1}$ and $\Delta z = \frac{1}{n_z-2}$. Furthermore, we take $n_x=N$ and $n_z=N+1$, so that $h=\Delta x = \Delta z$.
Then a second order finite difference discretization for~\eqref{eq:MM1aligned} and~\eqref{eq:MM2aligned} can be directly determined as
\begin{equation}\label{eq:MM1aligned_discret}
\left\{
\begin{array}{ll}
\begin{split}
 (-q_{i,j-1} + 2q_{i,j} - q_{i,j+1}) + \\
  (-\phi_{i-1,j} +  2\phi_{i,j}  - \phi_{i+1,j}) &
  = h^2 f_{i,j},
\end{split}
 & (i,j)\in[2,n_x-1]\times[2,n_z-1],\\[3mm]
q_{1,j}=q_{n_x,j}=0,& j\in[1,n_z],\\[3mm]
q_{i,1}=0,& i\in[2,n_x-1],\\[3mm]
-q_{i,n_z-1} + q_{i,n_z}=0,& i\in[2,n_x-1],
\end{array}
\right.
\end{equation}
and
\begin{equation}\label{eq:MM2aligned_discret}
\left\{
\begin{array}{ll}
\begin{split}
\varepsilon(q_{i,j-1} - 2q_{i,j} + q_{i,j+1}) + \\
(-\phi_{i,j-1} + 2\phi_{i,j} - &\phi_{i,j+1})   = 0,
\end{split}
& (i,j)\in[2,n_x-1]\times[2,n_z-1],\\[3mm]
\phi_{1,j}=\phi_{n_x,j}=0,& j\in[1,n_z],\\[3mm]
-\phi_{i,n_z-1} + \phi_{i,n_z}=0,& i\in[2,n_x-1],\\[3mm]
q_{i,1} - q_{i,2} =0,& i\in[2,n_x-1],
\end{array}
\right.
\end{equation}
where $\phi_{i,j}$ and $q_{i,j}$ are numerical approximations of $\phi(x_i,z_j)$ and $q(x_i,z_j)$ respectively, and $f_{i,j}= f(x_i,z_j)$.
Now, to construct the linear system to~\eqref{eq:MM1aligned_discret}-\eqref{eq:MM2aligned_discret}, we define a lexicographic, that is to count first in $z$ direction then in $x$ direction. Denote $iline$ the index of the linear system, we then define a function 'lexico' to convert from the couple $(i,j)$ to $iline$ as
\begin{equation}
iline = \text{lexico}(i,j) = j + (i-1)n_z.
\end{equation}
So the linear system corresponding to~\eqref{eq:MM1aligned_discret}-\eqref{eq:MM2aligned_discret} is given by
\begin{equation*}
\begin{pmatrix}
A_3 & B \\
A_2 & A_1
\end{pmatrix}
\begin{pmatrix}
Q \\
\Phi
 \end{pmatrix}
=
\begin{pmatrix}
F \\
0
\end{pmatrix}.
\end{equation*}
At last, in order to get a more symmetrical matrix forms, we make some algebraic operations, which yield in the first equation of~\eqref{eq:MM1aligned_discret} as
\begin{equation*}
\begin{array}{llll}
 -\phi_{i-1,j} + 2\phi_{i,j} - \phi_{i+1,j} & \to &  2\phi_{i,j} - \phi_{i+1,j},  & (i,j)\in\{2\}\times[2,n_z-1],\\[3mm]
-\phi_{i-1,j} + 2\phi_{i,j} - \phi_{i+1,j} & \to &   -\phi_{i-1,j} + 2\phi_{i,j}  ,& (i,j)\in\{n_x-1\}\times[2,n_z-1],\\[3mm]
 -q_{i,j-1} + 2q_{i,j} - q_{i,j+1} & \to & 2q_{i,j} - q_{i,j+1}, & (i,j)\in[2,n_x-1]\times \{2\}.
\end{array}
\end{equation*}
Thanks to these  algebraic operations, the matrices $A_3$ and $B$ become symmetric.
\subsection{Non-aligned case}\label{sec:discrete_nonaligned}

Now let us concentrate discretization for the general MMAP scheme~\eqref{eq:MM1}-\eqref{eq:MM2}. Again, this time we will consider a second order finite difference method.

Firstly, the mesh is defined different from the previous subsection, since we have to discretize cross derivative on all interior grids. For this, the mesh is defined by
\begin{equation*}
\left\{
\begin{array}{ll}
x_i = (i-1)\Delta x,\ i\in\{1,\dots,n_x\},\\[3mm]
z_j = (j-1)\Delta z,\ j\in\{1,\dots,n_z\},
\end{array}
\right.
\end{equation*}
where
$\Delta x = \frac{1}{n_x-1}$ and $\Delta z = \frac{1}{n_z-1}$. Again, we take $n_x=n_z=N$, so that $h=\Delta x = \Delta z$.

Secondly, we notice that
\begin{equation*}
\Delta \phi = \Delta_\bot\phi + \Delta_\|\phi,
\end{equation*}
its corresponding discretization on interior grids is
\begin{equation*}
\Delta^h \phi_{i,j} = \Delta^h_\bot\phi_{i,j} + \Delta^h_\|\phi_{i,j},\quad (i,j)\in\{2,\dots,n_x-1\}\times\{2,\dots,n_z-1\}.
\end{equation*}
Using the usual central discretization, $\Delta \phi$ can be approximated as
\begin{equation*}
h^2\Delta^h \phi_{i,j} = \phi_{i+1,j} +  \phi_{i-1,j} +  \phi_{i,j+1} +  \phi_{i,j-1} - 4  \phi_{i,j} .
\end{equation*}
Then, let us denote $\mathbf{b}\otimes\mathbf{b}$ by
\begin{equation*}
\mathbf{b}\otimes\mathbf{b} =
\begin{pmatrix}
b^{11} & b^{12} \\
b^{21} & b^{22}
\end{pmatrix},
\end{equation*}
thus we can recast $\Delta_\|\phi$ as
\begin{equation*}
\Delta_\|\phi = \partial_x G + \partial_z H,
\end{equation*}
where $G = b^{11}\partial_x\phi + b^{12}\partial_z\phi$, $H = b^{21}\partial_x\phi + b^{22}\partial_z\phi$.
Similarly, we get approximate $\Delta_\|\phi$ as
\begin{equation*}
h^2\Delta^h_\| \phi_{i,j} = (G^h_{i+1/2,j} - G^h_{i-1/2,j}) + (H^h_{i,j+1/2} - H^h_{i,j-1/2}),
\end{equation*}
where
\begin{eqnarray*}
G^h_{i+1/2,j} &=& b^{11}_{i+1/2,j}(\phi_{i+1,j} - \phi_{i,j}) + \frac{1}{4}b^{12}_{i+1/2,j}(\phi_{i+1,j+1} + \phi_{i,j+1} - \phi_{i+1,j-1} - \phi_{i,j-1}), \\
H^h_{i,j+1/2} &=&  \frac{1}{4}b^{21}_{i+1/2,j}(\phi_{i+1,j+1} + \phi_{i+1,j} - \phi_{i-1,j+1} - \phi_{i-1,j}) + b^{22}_{i,j+1/2}(\phi_{i,j+1} - \phi_{i,j}).
\end{eqnarray*}
At last, we obtain immediately approximated $\Delta_\bot\phi$ by
\begin{equation*}
 \Delta^h_\bot\phi_{i,j} = \Delta^h \phi_{i,j} - \Delta^h_\|\phi_{i,j}.
\end{equation*}
Thanks to these discrete operators, the sub-blocks $A_1$, $A_2$, $A_3$ and $B$ from the linear system~\eqref{eq:MMAPlinsyst} correspond to $-\Delta^h_\|\phi$, $\varepsilon\Delta^h_\|q$, $-\Delta^h_\|q$ and $-\Delta^h_\bot\phi$ respectively on interior grids.
Finally, let us consider discretization for boundary conditions.
From~\eqref{eq:MM1}-\eqref{eq:MM2}, we see that it is necessary to approximate $\mathbf{n}\cdot\nabla_\bot\phi$ and $\mathbf{n}\cdot\nabla_\|q$ on $[0,1]\times\{0,1\}$. We only present the approximation on boundary  $[0,1]\times\{0\}$, the one on boundary  $[0,1]\times\{1\}$ is similar.
Notice that the normal direction on  boundary  $[0,1]\times\{0\}$ is $\mathbf{n}=(0,1)^T$,  we thus has
\begin{eqnarray*}
\mathbf{n}\cdot\nabla_\bot\phi & = & b^{11} \partial_z \phi - b^{12}\partial_x \phi,\\
\mathbf{n}\cdot\nabla_\|q & = & b^{21} \partial_x q + b^{22}\partial_z q.
\end{eqnarray*}
Now applying center discretizations, we get a second approximation for above differential operations
\begin{eqnarray*}
h \mathbf{n}\cdot\nabla^h_\bot\phi_{i,1} & = & \frac{1}{2}b^{11}_{i,1}(-3\phi_{i,1} + 4\phi_{i,2} - \phi_{i,3}) - \frac{1}{2}b^{12}_{i,1}(\phi_{i+1,1} - \phi_{i-1,1}), \\
h \mathbf{n}\cdot\nabla^h_\| q_{i,1} & = & \frac{1}{2}b^{21}_{i,1}(q_{i+1,1} - q_{i-1,1}) + \frac{1}{2}b^{22}_{i,1}(-3q_{i,1} + 4q_{i,2} - q_{i,3}).
\end{eqnarray*}
At last, substituting  $\mathbf{n}\cdot\nabla_\bot\phi$ and $\mathbf{n}\cdot\nabla_\|q$ in~\eqref{eq:MM1}-\eqref{eq:MM2} by the above approximation, we obtain desired second order finite difference schemes.
\section{The property of the approximate Schur complement $S_4$}\label{sec:S4}
In this part, for aligned case we will show the approximate Schur complement $S_4$ of Section~\ref{sec:blockprecon} is precisely the exact one $E$.
We will consider the finite difference discretization described in Appendix~\ref{sec:discrete_aligned}, thus the relevant matrices are given as follows:

\begin{align*}
A_1 =
\begin{pmatrix}
I& \\
& a_1  \\
&& \ddots \\
&&&a_1\\
&&&& I
\end{pmatrix}
,
a_1 =
\begin{pmatrix}
0& \\
-1& 2 & -1  \\
&\ddots&\ddots& \ddots \\
&&-1&2&-1\\
&&&-1& 1
\end{pmatrix}.
\end{align*}

\begin{align*}
A_2 =
\begin{pmatrix}
0& \\
& a_2 \\
&& \ddots \\
&&&a_2\\
&&&& 0
\end{pmatrix}
,
a_2 =
\begin{pmatrix}
1& -1\\
\varepsilon& -2\varepsilon & \varepsilon  \\
&\ddots&\ddots& \ddots \\
&&\varepsilon&-2\varepsilon&\varepsilon\\
&&&& 0
\end{pmatrix}.
\end{align*}

\begin{align*}
A_3 =
\begin{pmatrix}
I& \\
& a_3  \\
&& \ddots \\
&&&a_3\\
&&&& I
\end{pmatrix}
,
a_3 =
\begin{pmatrix}
1& \\
& 2 & -1  \\
&-1& 2 & -1  \\
&&\ddots&\ddots& \ddots \\
&&&-1&2&-1\\
&&&&-1& 1
\end{pmatrix}.
\end{align*}

\begin{align*}
B =
\begin{pmatrix}
I& \\
-b_0& 2b_0& -b_0  \\
&& \ddots \\
&&-b_0&2b_0&-b_0\\
&&&& I
\end{pmatrix}
,
b_0 =
\begin{pmatrix}
0& \\
& 1 &   \\
&&\ddots&  \\
&&&1&\\
&&&& 0
\end{pmatrix}.
\end{align*}

Thanks to the results on tridiagonal matrices~\cite{CMD2007On}, we find
\begin{lemma}\label{lem:A3}
The inverse of the matrix $A_3$ is as follows
\begin{align*}
A_3^{-1}=\begin{pmatrix}
1& 0\\
0& C
\end{pmatrix},
C = \left[ \min\{i,j\} \right]_{i,j=1,\dots,N-1}.
\end{align*}
\end{lemma}

\begin{corollary}\label{cor:B}
From Lemma~\ref{lem:A3}, we immediately have
\begin{align*}
-A_2A_3^{-1}B =
\begin{pmatrix}
0& \\
-b_1& 2b_1& -b_1  \\
&& \ddots \\
&&-b_1&2b_1&-b_1\\
&&&& 0
\end{pmatrix}
,
b_1 =
\begin{pmatrix}
0& 1& \cdots&1&0\\
& \varepsilon&   \\
&&\ddots&  \\
&&&\varepsilon&\\
&&&& 0
\end{pmatrix}.
\end{align*}
\end{corollary}

Finally,  from Corollary~\ref{cor:B}, we deduce a consequence of the approximate Schur complement $S_4$ of Section~\ref{sec:blockprecon}.
\begin{proposition}
The approximate Schur complement $S_4$ is equal to the exact Schur complement $E$, {\it i.e.}
\begin{align*}
S_4 = A_1 - A_2A_3^{-1}B.
\end{align*}
\end{proposition}

We remark that if we substitute the first row of the matrix $b_1$ by null row, then the matrix $-A_2A_3^{-1}B$ becomes $\varepsilon B$. That is to say the only difference between the exact Schur complement $E$ and $A_1+\varepsilon B$ comes from these rows, which correspond to numerical approximation on the boundary $[0,1]\times\{0\}$. This fact for aligned case tells us that if we construct a better approximation of these rows,
we can obtain a better approximate Schur complement, and then the preconditioned Krylov solver converges more rapidly.

\end{document}